\documentclass[a4paper,10pt]{amsart}

\usepackage{amsfonts,amssymb}
\newtheorem{proposition}{Proposition}[section]
\newtheorem{lemma}[proposition]{Lemma}
\newtheorem{corollary}[proposition]{Corollary}
\newtheorem{theorem}[proposition]{Theorem}
\theoremstyle{definition}

\newtheorem{example}[proposition]{Example}

\theoremstyle{remark}
\newtheorem{remark}[proposition]{Remark}

\newcommand{\thlabel}[1]{\label{th:#1}}
\newcommand{\thref}[1]{Theorem~\ref{th:#1}}
\newcommand{\selabel}[1]{\label{se:#1}}
\newcommand{\seref}[1]{Section~\ref{se:#1}}
\newcommand{\lelabel}[1]{\label{le:#1}}
\newcommand{\leref}[1]{Lemma~\ref{le:#1}}
\newcommand{\prlabel}[1]{\label{pr:#1}}
\newcommand{\prref}[1]{Proposition~\ref{pr:#1}}

\newcommand{\eqlabel}[1]{\label{eq:#1}}
\newcommand{\equref}[1]{(\ref{eq:#1})}

\def\ra{\rightarrow}

\def\Id{{\rm Id}}

\def\ot{\otimes}
\def\va{\varepsilon}

\def\un{\underline}

\def\l{\lambda}

\def\va{\varepsilon}

\def\tr{\triangleright}

\def\ra{\rightarrow}

\def\cal{\mathcal}
\def\un{\underline}

\newcommand{\Cc}{\cal C}

\def\equal#1{\smash{\mathop{=}\limits^{#1}}}

\def\equalupdown#1#2{\smash{\mathop{=}\limits^{#1}\limits_{#2}}}

 \newcount\grcalca
 \newcount\grcalcb
 \newcount\grcalcc
 \newcount\grcalcd
 \newcount\grcalce
 \newcount\grcalcf
 \newcount\grcalcg
 \newcount\grcalch
 \newcount\grcalci
 \newcount\grrow
 \newcount\grcolumn
 \newcount\grwidth
 \newcount\factor
 \newcount\hfactor
 \newcount\qfactor
 \newcount\tfactor
 \newcount\sfactor
 \newcount\dfactor
 \hfactor = \factor
 \divide    \hfactor by 2
 \qfactor = \hfactor
 \divide    \qfactor by 2
 \tfactor = \factor
 \divide    \tfactor by 3
 \sfactor = \factor
 \divide    \sfactor by 6
 \dfactor = \factor
 \divide    \dfactor by 12

 \newcommand{\gbeg}[2]{
   \unitlength=1pt
   \grrow = #2
   \grcolumn = 0
   \grcalca = #1
   \grcalcb = #2
   \multiply \grcalca by \factor
   \grwidth = \grcalca
   \multiply \grcalcb by \factor
   \begin{minipage}{\grcalca pt}
   \begin{picture}(\grcalca,\grcalcb)
   \advance \grcalcb by -\factor
   \put(0, \grcalcb){\line(1,0){\grwidth}} }

 \newcommand{\gend}{
   \put(0, \factor){\line(1,0){\grwidth}}
   \end{picture}
   {\vskip2.5ex}
   \end{minipage} }

 \newcommand{\gnl}{
   \advance \grrow by -1
   \grcolumn = 0}

 \newcommand{\gvac}[1]{       
   \advance \grcolumn by #1} 
 
 \newcommand{\gcl}[1]{
   \grcalca = \grcolumn
   \multiply \grcalca by \factor
   \advance \grcalca by \hfactor
   \grcalcb = \grrow
   \multiply \grcalcb by \factor
   \grcalcc = #1
   \multiply \grcalcc by \factor
   \put(\grcalca,\grcalcb) {\line(0,-1){\grcalcc}} 
   \advance \grcolumn by 1}

 \newcommand{\gcn}[4]{
   \grcalca = \grcolumn
   \multiply \grcalca by \factor
   \grcalci = #3
   \multiply \grcalci by \hfactor
   \advance \grcalca by \grcalci
   \grcalcb = \grcolumn
   \multiply \grcalcb by \factor 
   \grcalci = #3
   \advance \grcalci by #4
   \multiply \grcalci by \qfactor
   \advance \grcalcb by \grcalci
   \grcalcc = \grcolumn
   \multiply \grcalcc by \factor
   \grcalci = #4
   \multiply \grcalci by \hfactor
   \advance \grcalcc by \grcalci
   \grcalcd = \grrow
   \multiply \grcalcd by \factor 
   \grcalce = \grrow
   \multiply \grcalce by \factor 
   \grcalci = #2
   \multiply \grcalci by \tfactor
   \advance \grcalce by -\grcalci
   \grcalcf = \grrow
   \multiply \grcalcf by \factor 
   \grcalci = #2
   \multiply \grcalci by \hfactor
   \advance \grcalcf by -\grcalci
   \grcalcg = \grrow
   \multiply \grcalcg by \factor 
   \grcalci = #2
   \multiply \grcalci by \tfactor
   \multiply \grcalci by 2
   \advance \grcalcg by -\grcalci
   \grcalch = \grrow
   \advance \grcalch by -#2
   \multiply \grcalch by \factor 
   \qbezier(\grcalca,\grcalcd)(\grcalca,\grcalce)(\grcalcb,\grcalcf) 
   \qbezier(\grcalcb,\grcalcf)(\grcalcc,\grcalcg)(\grcalcc,\grcalch) 
   \advance \grcolumn by #1}

 \newcommand{\gnot}[1]{
   \grcalca = \grcolumn
   \multiply \grcalca by \factor
   \advance \grcalca by \hfactor
   \grcalcb = \grrow
   \multiply \grcalcb by \factor
   \advance \grcalcb by -\hfactor
   \put(\grcalca,\grcalcb) {\makebox(0,0){$\scriptstyle #1$}} }

 \newcommand{\got}[2]{
   \grcalca = \grcolumn
   \multiply \grcalca by \factor
   \grcalcc = #1
   \multiply \grcalcc by \hfactor
   \advance \grcalca by \grcalcc
   \grcalcb = \grrow
   \multiply \grcalcb by \factor
   \advance \grcalcb by -\tfactor
   \advance \grcalcb by -\tfactor
   \put(\grcalca,\grcalcb){\makebox(0,0)[b]{$#2$}}
   \advance \grcolumn by #1}

 \newcommand{\gob}[2]{
   \grcalca = \grcolumn
   \multiply \grcalca by \factor
   \grcalcc = #1
   \multiply \grcalcc by \hfactor
   \advance \grcalca by \grcalcc
   \put(\grcalca,0){\makebox(0,0)[b]{$#2$}}
   \advance \grcolumn by #1}

 \newcommand{\gmu}{  
   \grcalca = \grcolumn
   \advance \grcalca by 1
   \multiply \grcalca by \factor
   \grcalcb = \grrow
   \multiply \grcalcb by \factor
   \grcalcc = \factor
   \advance \grcalcc by \hfactor
   \put(\grcalca,\grcalcb){\oval(\factor,\grcalcc)[b]}
   \advance \grcalcb by -\hfactor
   \advance \grcalcb by -\qfactor
   \put(\grcalca,\grcalcb) {\line(0,-1){\qfactor}} 
   \advance \grcolumn by 2}

 \newcommand{\gcmu}{   
   \grcalca = \grcolumn
   \advance \grcalca by 1
   \multiply \grcalca by \factor
   \grcalcb = \grrow
   \advance \grcalcb by -1
   \multiply \grcalcb by \factor
   \grcalcc = \factor
   \advance \grcalcc by \hfactor
   \put(\grcalca,\grcalcb){\oval(\factor,\grcalcc)[t]}
   \advance \grcalcb by \factor
   \put(\grcalca,\grcalcb) {\line(0,-1){\qfactor}} 
   \advance \grcolumn by 2}

 \newcommand{\glm}{
   \grcalca = \grcolumn
   \multiply \grcalca by \factor
   \advance \grcalca by \hfactor
   \grcalcb = \grcalca
   \advance \grcalcb by \factor
   \grcalcc = \grrow
   \multiply \grcalcc by \factor
   \grcalcd = \grcalcc
   \advance \grcalcd by -\tfactor
   \grcalce = \grcalcd
   \advance \grcalce by -\tfactor
   \put(\grcalca, \grcalcc){\line(0,-1){\tfactor}}
   \put(\grcalca, \grcalcd){\line(1,0){\factor}}
   \put(\grcalca, \grcalcd){\line(3,-1){\factor}}
   \put(\grcalcb, \grcalcc){\line(0,-1){\factor}}
   \advance \grcolumn by 2}

 \newcommand{\grm}{
   \grcalcb = \grcolumn
   \multiply \grcalcb by \factor
   \advance \grcalcb by \hfactor
   \grcalca = \grcalcb
   \advance \grcalca by \factor
   \grcalcc = \grrow
   \multiply \grcalcc by \factor
   \grcalcd = \grcalcc
   \advance \grcalcd by -\tfactor
   \grcalce = \grcalcd
   \advance \grcalce by -\tfactor
   \put(\grcalca, \grcalcc){\line(0,-1){\tfactor}}
   \put(\grcalca, \grcalcd){\line(-1,0){\factor}}
   \put(\grcalca, \grcalcd){\line(-3,-1){\factor}}
   \put(\grcalcb, \grcalcc){\line(0,-1){\factor}}
   \advance \grcolumn by 2}

 \newcommand{\glcm}{
   \grcalca = \grcolumn
   \multiply \grcalca by \factor
   \advance \grcalca by \hfactor
   \grcalcb = \grcalca
   \advance \grcalcb by \factor
   \grcalcc = \grrow
   \advance \grcalcc by -1
   \multiply \grcalcc by \factor
   \grcalcd = \grcalcc
   \advance \grcalcd by \tfactor
   \grcalce = \grcalcd
   \advance \grcalce by \tfactor
   \put(\grcalca, \grcalcc){\line(0,1){\tfactor}}
   \put(\grcalca, \grcalcd){\line(1,0){\factor}}
   \put(\grcalca, \grcalcd){\line(3,1){\factor}}
   \put(\grcalcb, \grcalcc){\line(0,1){\factor}}
   \advance \grcolumn by 2}

 \newcommand{\grcm}{
   \grcalcb = \grcolumn
   \multiply \grcalcb by \factor
   \advance \grcalcb by \hfactor
   \grcalca = \grcalcb
   \advance \grcalca by \factor
   \grcalcc = \grrow
   \advance \grcalcc by -1
   \multiply \grcalcc by \factor
   \grcalcd = \grcalcc
   \advance \grcalcd by \tfactor
   \grcalce = \grcalcd
   \advance \grcalce by \tfactor
   \put(\grcalca, \grcalcc){\line(0,1){\tfactor}}
   \put(\grcalca, \grcalcd){\line(-1,0){\factor}}
   \put(\grcalca, \grcalcd){\line(-3,1){\factor}}
   \put(\grcalcb, \grcalcc){\line(0,1){\factor}}
   \advance \grcolumn by 2}

 \newcommand{\gwmu}[1]{    
   \grcalca = \grcolumn
   \multiply \grcalca by \factor
   \grcalcd = \hfactor
   \multiply \grcalcd by #1
   \advance \grcalca by \grcalcd
   \grcalcb = \grrow
   \multiply \grcalcb by \factor
   \grcalcc = \factor
   \advance \grcalcc by \hfactor
   \grcalcd = #1
   \advance \grcalcd by -1
   \multiply \grcalcd by \factor
   \put(\grcalca,\grcalcb){\oval(\grcalcd,\grcalcc)[b]}
   \advance \grcalcb by -\hfactor
   \advance \grcalcb by -\qfactor
   \put(\grcalca,\grcalcb) {\line(0,-1){\qfactor}} 
   \advance \grcolumn by #1}

 \newcommand{\gwcm}[1]{   
   \grcalca = \grcolumn
   \multiply \grcalca by \factor
   \grcalcd = \hfactor
   \multiply \grcalcd by #1
   \advance \grcalca by \grcalcd
   \grcalcb = \grrow
   \advance \grcalcb by -1
   \multiply \grcalcb by \factor
   \grcalcc = \factor
   \advance \grcalcc by \hfactor
   \grcalcd = #1
   \advance \grcalcd by -1
   \multiply \grcalcd by \factor
   \put(\grcalca,\grcalcb){\oval(\grcalcd,\grcalcc)[t]}
   \advance \grcalcb by \factor
   \put(\grcalca,\grcalcb) {\line(0,-1){\qfactor}} 
   \advance \grcolumn by #1}

 \newcommand{\gwmuc}[1]{    
   \grcalca = \grcolumn
   \multiply \grcalca by \factor
   \advance \grcalca by \hfactor
   \grcalcb = \grrow
   \multiply \grcalcb by \factor
   \grcalcc = #1
   \advance \grcalcc by -1
   \multiply \grcalcc by \factor
   \put(\grcalca,\grcalcb){\line(1,0){\grcalcc}}
   \advance \grcalca by -\hfactor
   \grcalcd = \hfactor
   \multiply \grcalcd by #1
   \advance \grcalca by \grcalcd
   \grcalcc = \factor
   \advance \grcalcc by \hfactor
   \grcalcd = #1
   \advance \grcalcd by -1
   \multiply \grcalcd by \factor
   \put(\grcalca,\grcalcb){\oval(\grcalcd,\grcalcc)[b]}
   \advance \grcalcb by -\hfactor
   \advance \grcalcb by -\qfactor
   \put(\grcalca,\grcalcb) {\line(0,-1){\qfactor}} 
   \advance \grcolumn by #1}

 \newcommand{\gwcmc}[1]{   
   \grcalca = \grcolumn
   \multiply \grcalca by \factor
   \advance \grcalca by \hfactor
   \grcalcb = \grrow
   \multiply \grcalcb by \factor
   \advance \grcalcb by -\factor
   \grcalcc = #1
   \advance \grcalcc by -1
   \multiply \grcalcc by \factor
   \put(\grcalca,\grcalcb){\line(1,0){\grcalcc}}
   \grcalcd = #1
   \advance \grcalcd by -1
   \multiply \grcalcd by \hfactor
   \advance \grcalca by \grcalcd
   \grcalcc = \factor
   \advance \grcalcc by \hfactor
   \grcalcd = #1
   \advance \grcalcd by -1
   \multiply \grcalcd by \factor
   \put(\grcalca,\grcalcb){\oval(\grcalcd,\grcalcc)[t]}
   \advance \grcalcb by \factor
   \put(\grcalca,\grcalcb) {\line(0,-1){\qfactor}} 
   \advance \grcolumn by #1}

 \newcommand{\gev}{  
   \grcalca = \grcolumn
   \advance \grcalca by 1
   \multiply \grcalca by \factor
   \grcalcb = \grrow
   \multiply \grcalcb by \factor
   \grcalcc = \factor
   \advance \grcalcc by \hfactor
   \put(\grcalca,\grcalcb){\oval(\factor,\grcalcc)[b]}
   \advance \grcolumn by 2}

 \newcommand{\gdb}{   
   \grcalca = \grcolumn
   \advance \grcalca by 1
   \multiply \grcalca by \factor
   \grcalcb = \grrow
   \advance \grcalcb by -1
   \multiply \grcalcb by \factor
   \grcalcc = \factor
   \advance \grcalcc by \hfactor
   \put(\grcalca,\grcalcb){\oval(\factor,\grcalcc)[t]}
   \advance \grcolumn by 2}

 \newcommand{\gwev}[1]{    
   \grcalca = \grcolumn
   \multiply \grcalca by \factor
   \grcalcd = \hfactor
   \multiply \grcalcd by #1
   \advance \grcalca by \grcalcd
   \grcalcb = \grrow
   \multiply \grcalcb by \factor
   \grcalcc = \factor
   \advance \grcalcc by \hfactor
   \grcalcd = #1
   \advance \grcalcd by -1
   \multiply \grcalcd by \factor
   \put(\grcalca,\grcalcb){\oval(\grcalcd,\grcalcc)[b]}
   \advance \grcolumn by #1}

 \newcommand{\gwdb}[1]{   
   \grcalca = \grcolumn
   \multiply \grcalca by \factor
   \grcalcd = \hfactor
   \multiply \grcalcd by #1
   \advance \grcalca by \grcalcd
   \grcalcb = \grrow
   \advance \grcalcb by -1
   \multiply \grcalcb by \factor
   \grcalcc = \factor
   \advance \grcalcc by \hfactor
   \grcalcd = #1
   \advance \grcalcd by -1
   \multiply \grcalcd by \factor
   \put(\grcalca,\grcalcb){\oval(\grcalcd,\grcalcc)[t]}
   \advance \grcolumn by #1}

 \newcommand{\gbr}{
   \grcalca = \grcolumn
   \multiply \grcalca by \factor
   \advance \grcalca by \hfactor
   \grcalcb = \grcalca
   \advance \grcalcb by \hfactor
   \grcalcc = \grcalca
   \advance \grcalcc by \factor
   \grcalcd = \grrow
   \multiply \grcalcd by \factor
   \grcalce = \grcalcd
   \advance \grcalce by -\tfactor
   \grcalcf = \grcalcd
   \advance \grcalcf by -\hfactor
   \grcalcg = \grcalce
   \advance \grcalcg by -\tfactor
   \grcalch = \grcalcd
   \advance \grcalch by -\factor
   \qbezier(\grcalca,\grcalcd)(\grcalca,\grcalce)(\grcalcb,\grcalcf) 
   \qbezier(\grcalcb,\grcalcf)(\grcalcc,\grcalcg)(\grcalcc,\grcalch) 
   \advance \grcalcf by -\dfactor
   \advance \grcalcb by -\sfactor
   \qbezier(\grcalca,\grcalch)(\grcalca,\grcalcg)(\grcalcb,\grcalcf) 
   \advance \grcalcf by \sfactor
   \advance \grcalcb by \tfactor
   \qbezier(\grcalcc,\grcalcd)(\grcalcc,\grcalce)(\grcalcb,\grcalcf) 
   \advance \grcolumn by 2}

 \newcommand{\gibr}{
   \grcalca = \grcolumn
   \multiply \grcalca by \factor
   \advance \grcalca by \hfactor
   \grcalcb = \grcalca
   \advance \grcalcb by \hfactor
   \grcalcc = \grcalca
   \advance \grcalcc by \factor
   \grcalcd = \grrow
   \multiply \grcalcd by \factor
   \grcalce = \grcalcd
   \advance \grcalce by -\tfactor
   \grcalcf = \grcalcd
   \advance \grcalcf by -\hfactor
   \grcalcg = \grcalce
   \advance \grcalcg by -\tfactor
   \grcalch = \grcalcd
   \advance \grcalch by -\factor
   \qbezier(\grcalcc,\grcalcd)(\grcalcc,\grcalce)(\grcalcb,\grcalcf) 
   \qbezier(\grcalcb,\grcalcf)(\grcalca,\grcalcg)(\grcalca,\grcalch) 
   \advance \grcalcf by -\dfactor
   \advance \grcalcb by \sfactor
   \qbezier(\grcalcc,\grcalch)(\grcalcc,\grcalcg)(\grcalcb,\grcalcf) 
   \advance \grcalcf by \sfactor
   \advance \grcalcb by -\tfactor
   \qbezier(\grcalca,\grcalcd)(\grcalca,\grcalce)(\grcalcb,\grcalcf) 
   \advance \grcolumn by 2}

 \newcommand{\gbrc}{
   \grcalca = \grcolumn
   \multiply \grcalca by \factor
   \advance \grcalca by \hfactor
   \grcalcb = \grcalca
   \advance \grcalcb by \hfactor
   \grcalcc = \grcalca
   \advance \grcalcc by \factor
   \grcalcd = \grrow
   \multiply \grcalcd by \factor
   \grcalce = \grcalcd
   \advance \grcalce by -\tfactor
   \grcalcf = \grcalcd
   \advance \grcalcf by -\hfactor
   \grcalcg = \grcalce
   \advance \grcalcg by -\tfactor
   \grcalch = \grcalcd
   \advance \grcalch by -\factor
   \put(\grcalcb,\grcalcf){\circle{\hfactor}}
   \qbezier(\grcalca,\grcalcd)(\grcalca,\grcalce)(\grcalcb,\grcalcf) 
   \qbezier(\grcalcb,\grcalcf)(\grcalcc,\grcalcg)(\grcalcc,\grcalch) 
   \advance \grcalcf by -\dfactor
   \advance \grcalcb by -\sfactor
   \qbezier(\grcalca,\grcalch)(\grcalca,\grcalcg)(\grcalcb,\grcalcf) 
   \advance \grcalcf by \sfactor
   \advance \grcalcb by \tfactor
   \qbezier(\grcalcc,\grcalcd)(\grcalcc,\grcalce)(\grcalcb,\grcalcf) 
   \advance \grcolumn by 2}

 \newcommand{\gibrc}{
   \grcalca = \grcolumn
   \multiply \grcalca by \factor
   \advance \grcalca by \hfactor
   \grcalcb = \grcalca
   \advance \grcalcb by \hfactor
   \grcalcc = \grcalca
   \advance \grcalcc by \factor
   \grcalcd = \grrow
   \multiply \grcalcd by \factor
   \grcalce = \grcalcd
   \advance \grcalce by -\tfactor
   \grcalcf = \grcalcd
   \advance \grcalcf by -\hfactor
   \grcalcg = \grcalce
   \advance \grcalcg by -\tfactor
   \grcalch = \grcalcd
   \advance \grcalch by -\factor
   \put(\grcalcb,\grcalcf){\circle{\hfactor}}
   \qbezier(\grcalcc,\grcalcd)(\grcalcc,\grcalce)(\grcalcb,\grcalcf) 
   \qbezier(\grcalcb,\grcalcf)(\grcalca,\grcalcg)(\grcalca,\grcalch) 
   \advance \grcalcf by -\dfactor
   \advance \grcalcb by \sfactor
   \qbezier(\grcalcc,\grcalch)(\grcalcc,\grcalcg)(\grcalcb,\grcalcf) 
   \advance \grcalcf by \sfactor
   \advance \grcalcb by -\tfactor
   \qbezier(\grcalca,\grcalcd)(\grcalca,\grcalce)(\grcalcb,\grcalcf) 
   \advance \grcolumn by 2} 

 \newcommand{\gu}[1]{
   \grcalca = \grcolumn
   \multiply \grcalca by \factor
   \grcalcd = \hfactor
   \multiply \grcalcd by #1
   \advance \grcalca by \grcalcd
   \grcalcb = \grrow
   \advance \grcalcb by -1
   \multiply \grcalcb by \factor
   \put(\grcalca,\grcalcb) {\line(0,1){\hfactor}} 
   \advance \grcalcb by \hfactor
   \put(\grcalca,\grcalcb) {\circle*{3}}
   \advance \grcolumn by #1}

 \newcommand{\gcu}[1]{
   \grcalca = \grcolumn
   \multiply \grcalca by \factor
   \grcalcd = \hfactor
   \multiply \grcalcd by #1
   \advance \grcalca by \grcalcd
   \grcalcb = \grrow
   \multiply \grcalcb by \factor
   \put(\grcalca,\grcalcb) {\line(0,-1){\hfactor}} 
   \advance \grcalcb by -\hfactor
   \put(\grcalca,\grcalcb) {\circle*{3}}
   \advance \grcolumn by #1}

 \newcommand{\gmp}[1]{
   \grcalca = \grcolumn
   \multiply \grcalca by \factor
   \advance \grcalca by \hfactor
   \grcalcb = \grrow
   \multiply \grcalcb by \factor
   \put(\grcalca,\grcalcb) {\line(0,-1){\dfactor}} 
   \advance \grcalcb by -\factor
   \put(\grcalca,\grcalcb) {\line(0,1){\dfactor}} 
   \advance \grcalcb by \hfactor
   \grcalcc = \factor
   \advance \grcalcc by -\qfactor
   \put(\grcalca,\grcalcb) {\circle{\grcalcc}}
   \put(\grcalca,\grcalcb) {\makebox(0,0){$\scriptstyle #1$}}
   \advance \grcolumn by 1}

 \newcommand{\gbmp}[1]{
   \grcalca = \grcolumn
   \multiply \grcalca by \factor
   \advance \grcalca by \hfactor
   \grcalcb = \grrow
   \multiply \grcalcb by \factor
   \put(\grcalca,\grcalcb) {\line(0,-1){\dfactor}} 
   \advance \grcalcb by -\factor
   \put(\grcalca,\grcalcb) {\line(0,1){\dfactor}} 
   \advance \grcalca by -\hfactor
   \advance \grcalca by \dfactor
   \advance \grcalcb by \dfactor
   \grcalcc = \factor
   \advance \grcalcc by -\sfactor
   \put(\grcalca,\grcalcb) {\framebox(\grcalcc,\grcalcc){$\scriptstyle #1$}}
   \advance \grcolumn by 1}

 \newcommand{\gbmpt}[1]{
   \grcalca = \grcolumn
   \multiply \grcalca by \factor
   \advance \grcalca by \hfactor
   \grcalcb = \grrow
   \multiply \grcalcb by \factor
   \put(\grcalca,\grcalcb) {\line(0,-1){\dfactor}} 
   \advance \grcalcb by -\factor
   \advance \grcalca by -\hfactor
   \advance \grcalca by \dfactor
   \advance \grcalcb by \dfactor
   \grcalcc = \factor
   \advance \grcalcc by -\sfactor
   \put(\grcalca,\grcalcb) {\framebox(\grcalcc,\grcalcc){$\scriptstyle #1$}}
   \advance \grcolumn by 1}

 \newcommand{\gbmpb}[1]{
   \grcalca = \grcolumn
   \multiply \grcalca by \factor
   \advance \grcalca by \hfactor
   \grcalcb = \grrow
   \multiply \grcalcb by \factor
   \advance \grcalcb by -\factor
   \put(\grcalca,\grcalcb) {\line(0,1){\dfactor}} 
   \advance \grcalca by -\hfactor
   \advance \grcalca by \dfactor
   \advance \grcalcb by \dfactor
   \grcalcc = \factor
   \advance \grcalcc by -\sfactor
   \put(\grcalca,\grcalcb) {\framebox(\grcalcc,\grcalcc){$\scriptstyle #1$}}
   \advance \grcolumn by 1}

 \newcommand{\gbmpn}[1]{
   \grcalca = \grcolumn
   \multiply \grcalca by \factor
   \advance \grcalca by \hfactor
   \grcalcb = \grrow
   \multiply \grcalcb by \factor
   \advance \grcalcb by -\factor
   \advance \grcalca by -\hfactor
   \advance \grcalca by \dfactor
   \advance \grcalcb by \dfactor
   \grcalcc = \factor
   \advance \grcalcc by -\sfactor
   \put(\grcalca,\grcalcb) {\framebox(\grcalcc,\grcalcc){$\scriptstyle #1$}}
   \advance \grcolumn by 1}

 \newcommand{\glmptb}{    
   \grcalca = \grcolumn
   \multiply \grcalca by \factor
   \advance \grcalca by \hfactor
   \grcalcb = \grrow
   \multiply \grcalcb by \factor
   \put(\grcalca,\grcalcb) {\line(0,-1){\dfactor}} 
   \advance \grcalcb by -\factor
   \put(\grcalca,\grcalcb) {\line(0,1){\dfactor}} 
   \advance \grcalca by -\hfactor
   \advance \grcalca by \dfactor
   \advance \grcalcb by \dfactor
   \put(\grcalca,\grcalcb) {\line(1,0){\factor}} 
   \advance \grcalcb by \factor
   \advance \grcalcb by -\sfactor
   \put(\grcalca,\grcalcb) {\line(1,0){\factor}} 
   \grcalcc = \factor
   \advance \grcalcc by -\sfactor
   \put(\grcalca,\grcalcb) {\line(0,-1){\grcalcc}} 
   \advance \grcolumn by 1}

 \newcommand{\glmpt}{    
   \grcalca = \grcolumn
   \multiply \grcalca by \factor
   \advance \grcalca by \hfactor
   \grcalcb = \grrow
   \multiply \grcalcb by \factor
   \put(\grcalca,\grcalcb) {\line(0,-1){\dfactor}} 
   \advance \grcalca by -\hfactor
   \advance \grcalca by \dfactor
   \advance \grcalcb by -\dfactor
   \put(\grcalca,\grcalcb) {\line(1,0){\factor}} 
   \advance \grcalcb by -\factor
   \advance \grcalcb by \sfactor
   \put(\grcalca,\grcalcb) {\line(1,0){\factor}} 
   \grcalcc = \factor
   \advance \grcalcc by -\sfactor
   \put(\grcalca,\grcalcb) {\line(0,1){\grcalcc}} 
   \advance \grcolumn by 1}

 \newcommand{\glmpb}{    
   \grcalca = \grcolumn
   \multiply \grcalca by \factor
   \advance \grcalca by \hfactor
   \grcalcb = \grrow
   \multiply \grcalcb by \factor
   \advance \grcalcb by -\factor
   \put(\grcalca,\grcalcb) {\line(0,1){\dfactor}} 
   \advance \grcalca by -\hfactor
   \advance \grcalca by \dfactor
   \advance \grcalcb by \dfactor
   \put(\grcalca,\grcalcb) {\line(1,0){\factor}} 
   \advance \grcalcb by \factor
   \advance \grcalcb by -\sfactor
   \put(\grcalca,\grcalcb) {\line(1,0){\factor}} 
   \grcalcc = \factor
   \advance \grcalcc by -\sfactor
   \put(\grcalca,\grcalcb) {\line(0,-1){\grcalcc}} 
   \advance \grcolumn by 1}

 \newcommand{\glmp}{    
   \grcalca = \grcolumn
   \multiply \grcalca by \factor
   \advance \grcalca by \dfactor
   \grcalcb = \grrow
   \multiply \grcalcb by \factor
   \advance \grcalcb by -\dfactor
   \put(\grcalca,\grcalcb) {\line(1,0){\factor}} 
   \advance \grcalcb by -\factor
   \advance \grcalcb by \sfactor
   \put(\grcalca,\grcalcb) {\line(1,0){\factor}} 
   \grcalcc = \factor
   \advance \grcalcc by -\sfactor
   \put(\grcalca,\grcalcb) {\line(0,1){\grcalcc}} 
   \advance \grcolumn by 1}

 \newcommand{\gcmptb}{    
   \grcalca = \grcolumn
   \multiply \grcalca by \factor
   \advance \grcalca by \hfactor
   \grcalcb = \grrow
   \multiply \grcalcb by \factor
   \put(\grcalca,\grcalcb) {\line(0,-1){\dfactor}} 
   \advance \grcalcb by -\factor
   \put(\grcalca,\grcalcb) {\line(0,1){\dfactor}} 
   \advance \grcalca by -\hfactor
   \advance \grcalcb by \dfactor
   \put(\grcalca,\grcalcb) {\line(1,0){\factor}} 
   \advance \grcalcb by \factor
   \advance \grcalcb by -\sfactor
   \put(\grcalca,\grcalcb) {\line(1,0){\factor}} 
   \advance \grcolumn by 1}

 \newcommand{\gcmpt}{    
   \grcalca = \grcolumn
   \multiply \grcalca by \factor
   \advance \grcalca by \hfactor
   \grcalcb = \grrow
   \multiply \grcalcb by \factor
   \put(\grcalca,\grcalcb) {\line(0,-1){\dfactor}} 
   \advance \grcalcb by -\factor
   \advance \grcalca by -\hfactor
   \advance \grcalcb by \dfactor
   \put(\grcalca,\grcalcb) {\line(1,0){\factor}} 
   \advance \grcalcb by \factor
   \advance \grcalcb by -\sfactor
   \put(\grcalca,\grcalcb) {\line(1,0){\factor}} 
   \advance \grcolumn by 1}

 \newcommand{\gcmpb}{    
   \grcalca = \grcolumn
   \multiply \grcalca by \factor
   \advance \grcalca by \hfactor
   \grcalcb = \grrow
   \multiply \grcalcb by \factor
   \advance \grcalcb by -\factor
   \put(\grcalca,\grcalcb) {\line(0,1){\dfactor}} 
   \advance \grcalca by -\hfactor
   \advance \grcalcb by \dfactor
   \put(\grcalca,\grcalcb) {\line(1,0){\factor}} 
   \advance \grcalcb by \factor
   \advance \grcalcb by -\sfactor
   \put(\grcalca,\grcalcb) {\line(1,0){\factor}} 
   \advance \grcolumn by 1}

 \newcommand{\gcmp}{    
   \grcalca = \grcolumn
   \multiply \grcalca by \factor
   \grcalcb = \grrow
   \multiply \grcalcb by \factor
   \advance \grcalcb by -\factor
   \advance \grcalcb by \dfactor
   \put(\grcalca,\grcalcb) {\line(1,0){\factor}} 
   \advance \grcalcb by \factor
   \advance \grcalcb by -\sfactor
   \put(\grcalca,\grcalcb) {\line(1,0){\factor}} 
   \advance \grcolumn by 1}

 \newcommand{\grmptb}{    
   \grcalca = \grcolumn
   \multiply \grcalca by \factor
   \advance \grcalca by \hfactor
   \grcalcb = \grrow
   \multiply \grcalcb by \factor
   \put(\grcalca,\grcalcb) {\line(0,-1){\dfactor}} 
   \advance \grcalcb by -\factor
   \put(\grcalca,\grcalcb) {\line(0,1){\dfactor}} 
   \advance \grcalca by \hfactor
   \advance \grcalca by -\dfactor
   \advance \grcalcb by \dfactor
   \put(\grcalca,\grcalcb) {\line(-1,0){\factor}} 
   \advance \grcalcb by \factor
   \advance \grcalcb by -\sfactor
   \put(\grcalca,\grcalcb) {\line(-1,0){\factor}} 
   \grcalcc = \factor
   \advance \grcalcc by -\sfactor
   \put(\grcalca,\grcalcb) {\line(0,-1){\grcalcc}} 
   \advance \grcolumn by 1}

 \newcommand{\grmpt}{    
   \grcalca = \grcolumn
   \multiply \grcalca by \factor
   \advance \grcalca by \hfactor
   \grcalcb = \grrow
   \multiply \grcalcb by \factor
   \put(\grcalca,\grcalcb) {\line(0,-1){\dfactor}} 
   \advance \grcalca by \hfactor
   \advance \grcalca by -\dfactor
   \advance \grcalcb by -\dfactor
   \put(\grcalca,\grcalcb) {\line(-1,0){\factor}} 
   \advance \grcalcb by -\factor
   \advance \grcalcb by \sfactor
   \put(\grcalca,\grcalcb) {\line(-1,0){\factor}} 
   \grcalcc = \factor
   \advance \grcalcc by -\sfactor
   \put(\grcalca,\grcalcb) {\line(0,1){\grcalcc}} 
   \advance \grcolumn by 1}

 \newcommand{\grmpb}{    
   \grcalca = \grcolumn
   \multiply \grcalca by \factor
   \advance \grcalca by \hfactor
   \grcalcb = \grrow
   \multiply \grcalcb by \factor
   \advance \grcalcb by -\factor
   \put(\grcalca,\grcalcb) {\line(0,1){\dfactor}} 
   \advance \grcalca by \hfactor
   \advance \grcalca by -\dfactor
   \advance \grcalcb by \dfactor
   \put(\grcalca,\grcalcb) {\line(-1,0){\factor}} 
   \advance \grcalcb by \factor
   \advance \grcalcb by -\sfactor
   \put(\grcalca,\grcalcb) {\line(-1,0){\factor}} 
   \grcalcc = \factor
   \advance \grcalcc by -\sfactor
   \put(\grcalca,\grcalcb) {\line(0,-1){\grcalcc}} 
   \advance \grcolumn by 1}

 \newcommand{\grmp}{    
   \grcalca = \grcolumn
   \multiply \grcalca by \factor
   \advance \grcalca by \factor
   \advance \grcalca by -\dfactor
   \grcalcb = \grrow
   \multiply \grcalcb by \factor
   \advance \grcalcb by -\dfactor
   \put(\grcalca,\grcalcb) {\line(-1,0){\factor}} 
   \advance \grcalcb by -\factor
   \advance \grcalcb by \sfactor
   \put(\grcalca,\grcalcb) {\line(-1,0){\factor}} 
   \grcalcc = \factor
   \advance \grcalcc by -\sfactor
   \put(\grcalca,\grcalcb) {\line(0,1){\grcalcc}} 
   \advance \grcolumn by 1}
\newcommand{\gsy}{
   \grcalca = \grcolumn
   \multiply \grcalca by \factor
   \advance \grcalca by \hfactor
   \grcalcb = \grcalca
   \advance \grcalcb by \hfactor
   \grcalcc = \grcalca
   \advance \grcalcc by \factor
   \grcalcd = \grrow
   \multiply \grcalcd by \factor
   \grcalce = \grcalcd
   \advance \grcalce by -\tfactor
   \grcalcf = \grcalcd
   \advance \grcalcf by -\hfactor
   \grcalcg = \grcalce
   \advance \grcalcg by -\tfactor
   \grcalch = \grcalcd
   \advance \grcalch by -\factor
   \qbezier(\grcalcc,\grcalcd)(\grcalcc,\grcalce)(\grcalcb,\grcalcf) 
   \qbezier(\grcalcb,\grcalcf)(\grcalca,\grcalcg)(\grcalca,\grcalch) 
   \advance \grcalcf by -\dfactor
   \advance \grcalcb by \sfactor
   \qbezier(\grcalcc,\grcalch)(\grcalcc,\grcalcg)(\grcalcb,\grcalcf) 
   \qbezier(\grcalca,\grcalcd)(\grcalca,\grcalce)(\grcalcb,\grcalcf) 
   \advance \grcolumn by 2}

 \newcommand{\gwmuh}[3]{    
   \grcalca = \grcolumn
   \multiply \grcalca by \factor
   \grcalcb = #2
   \advance \grcalcb by #3
   \multiply \grcalcb by \qfactor
   \advance \grcalca by \grcalcb
   \grcalcb = \grrow
   \multiply \grcalcb by \factor
   \grcalcc = #3
   \advance \grcalcc by -#2
   \multiply \grcalcc by \hfactor
   \grcalcd = \factor
   \advance \grcalcd by \hfactor
   \put(\grcalca,\grcalcb){\oval(\grcalcc,\grcalcd)[b]}
   \grcalca = \grcolumn
   \multiply \grcalca by \factor
   \grcalcc = #1
   \multiply \grcalcc by \hfactor
   \advance \grcalca by \grcalcc
   \advance \grcalcb by -\hfactor
   \advance \grcalcb by -\qfactor
   \put(\grcalca,\grcalcb) {\line(0,-1){\qfactor}} 
   \advance \grcolumn by #1}

 \newcommand{\gwcmh}[3]{   
   \grcalca = \grcolumn
   \multiply \grcalca by \factor
   \grcalcb = #2
   \advance \grcalcb by #3
   \multiply \grcalcb by \qfactor
   \advance \grcalca by \grcalcb
   \grcalcb = \grrow
   \advance \grcalcb by -1
   \multiply \grcalcb by \factor
   \grcalcc = #3
   \advance \grcalcc by -#2
   \multiply \grcalcc by \hfactor
   \grcalcd = \factor
   \advance \grcalcd by \hfactor
   \put(\grcalca,\grcalcb){\oval(\grcalcc,\grcalcd)[t]}
   \grcalca = \grcolumn
   \multiply \grcalca by \factor
   \grcalcc = #1
   \multiply \grcalcc by \hfactor
   \advance \grcalca by \grcalcc
   \advance \grcalcb by \factor
   \put(\grcalca,\grcalcb) {\line(0,-1){\qfactor}} 
   \advance \grcolumn by #1}

 \newcommand{\gsbox}[1]{
   \grcalca = \grcolumn
   \multiply \grcalca by \factor
   \grcalcb = \grrow
   \multiply \grcalcb by \factor
   \advance \grcalcb by -\factor
   \grcalcc = #1
   \multiply \grcalcc by \factor
   \grcalcd = \factor
   \put(\grcalca,\grcalcb){\framebox(\grcalcc,\grcalcd){}}}

\allowdisplaybreaks[4]

\begin{document}
\title[A monoidal structure]
{A monoidal structure on the category of relative Hopf modules}
\author{D. Bulacu}
\address{Faculty of Mathematics and Informatics, University
of Bucharest, Strada Academiei 14, RO-010014 Bucharest 1, Romania and
Faculty of Engineering, 
Vrije Universiteit Brussel, B-1050 Brussels, Belgium}
\email{daniel.bulacu@fmi.unibuc.ro}
\author{S. Caenepeel}
\address{Faculty of Engineering, 
Vrije Universiteit Brussel, B-1050 Brussels, Belgium}
\email{scaenepe@vub.ac.be}
\urladdr{http://homepages.vub.ac.be/\~{}scaenepe/}
\thanks{
The first author was supported by FWO-Vlaanderen
(FWO GP.045.09N), and by CNCSIS $479/2009$, code ID 1904.
This research is part of the FWO project G.0117.10 ``Equivariant Brauer groups
and Galois deformations". The authors also thank Bodo Pareigis for sharing his
``diagrams'' program.}
\subjclass[2010]{16T05, 18D10}

\keywords{Monoidal category, Relative Hopf module, Yetter-Drinfeld module, braided bialgebra}

\begin{abstract}
Let $B$ be a bialgebra, and $A$ a left $B$-comodule algebra in a braided monoidal
category $\Cc$, and assume that $A$ is also a coalgebra, with a not-necessarily
associative or unital left $B$-action. Then we can define a right $A$-action on the tensor
product of two relative Hopf modules, and this defines a monoidal structure on the category 
of relative Hopf modules if and only if $A$ is a bialgebra in the category of
left Yetter-Drinfeld modules over $B$. Some examples are given. 
\end{abstract}

\maketitle
\section*{Introduction}\selabel{intro}
It is well-known that the category  of corepresentations over a bialgebra $B$ in a braided
monoidal category $\Cc$ is monoidal. Now let $A$ be a left $B$-comodule algebra,
and consider the category of relative Hopf modules ${}^B\Cc_A$. A relative Hopf module
is always a left $B$-comodule, in fact we have a forgetful functor ${}^B\Cc_A\to {}^B\Cc$.
The following natural question arises: is there a monoidal structure on ${}^B\Cc_A$
that is compatible with the one on ${}^B\Cc$, by which we mean that the forgetful
functor is strongly monoidal.\\
Monoidal structures on a more general category, the category of Doi-Hopf modules, have been discussed 
in \cite{cob}, in the particular case where $\Cc$
is the category of vector spaces over a field $k$. A monoidal structure on ${}^B\Cc_A$
can be constructed if $A$ is a bialgebra and two additional compatibility conditions 
are satisfied. The aim of this paper is to present a more general result in the following
direction: we will no longer assume that $A$ is a bialgebra: it will be sufficient that
$A$ is at the same time an algebra and a coalgebra, and, moreover, we will assume
that we have a left $B$-action $B\ot A\ra A$, which is not assumed to be associative or unital 
from the beginning. These additional structures (coalgebra and $B$-action) on $A$
allow us to define a right $A$-action on the tensor product of two relative Hopf modules, 
and on the unit object $\un{1}$ of $\Cc$. Left $B$-coactions on these objects are
supplied using the monoidal structure of the category of left $B$-comodules.
Our main result, \thref{2.1}, states that the category of relative Hopf modules
with this additional structure, is a monoidal category if and only if
$A$ is a braided bialgebra, this is a bialgebra 
in the prebraided monoidal category ${}^B_B{\cal YD}$ of left Yetter-Drinfeld 
modules. In the case where $\Cc$ is the category of vector spaces, braided bialgebras
and Hopf algebras appeared in the theory of Hopf algebras with a projection \cite{rad}.
Observe also that braided Hopf algebras play an important role in the classification theory
of pointed Hopf algebras, see for example \cite{as1} for a survey and \cite{as} for the most
recent developments.\\
In \seref{3}, we discuss some particular situations and examples.
We first consider a quasitriangular Hopf algebra and its associated
enveloping algebra braided group $\un{H}$ as in \cite{maj}. It is well-known that
$\un{H}$ is a braided bialgebra, and therefore ${}^H{\cal M}_{\un{H}}$ is a monoidal category.
As a consequence, we find that the category of Long $H$-dimodules over a cocommutative
Hopf algebra is monoidal. A second example is provided by a coquasitriangular Hopf
algebra $(H, \sigma)$ and its associated (left) function algebra braided group $\un{H}$, which
is a braided bialgebra, so that the category of relative Hopf modules
${}^H{\cal M}_{\un{H}}$ is monoidal. In the case where $H$ is commutative and $\sigma$
is trivial, this category is identified with the category of Yetter-Drinfeld modules over $H$.\\
Finally, we look at the situation where the given left $B$-action on $A$ is trivial, and
we show that $(B,A,B)$ is a monoidal Doi-Hopf datum in the sense of \cite{cob}
if and only if $A$ is a braided bialgebra with trivial $B$-action.

\section{Preliminary results}\selabel{prelimres}
\setcounter{equation}{0}
\subsection{Braided monoidal categories}\selabel{1.1}
We assume that the reader is familiar with the basic theory of braided monoidal categories;
for details, we refer to \cite{kas, maj}.
In the sequel $\Cc$ is a (pre)braided monoidal category 
with tensor product $\ot : \Cc\times \Cc\ra \Cc$, 
unit object $\un{1}$ and (pre)braiding $c: \ot\ra \ot\circ \tau$, where 
$\tau: \Cc\times \Cc\ra \Cc\times \Cc$ is the twist functor. 
For any two objects $X, Y$ of $\Cc$ 
we denote $c_{X, Y}$ by 
$
\gbeg{2}{3}
\got{1}{X}\got{1}{Y}\gnl
\gbr\gnl
\gob{1}{Y}\gob{1}{X}
\gend
$.  
Recall that a (pre)braiding $c$ satisfies  
\begin{equation}\label{braiding}
c_{X, Y\ot Z}=
\gbeg{3}{4}
\got{1}{X}\got{1}{Y}\got{1}{Z}\gnl
\gbr\gcl{1}\gnl
\gcl{1}\gbr\gnl
\gob{1}{Y}\gob{1}{Z}\gob{1}{X}
\gend
\mbox{\hspace{4mm}and\hspace{4mm}}
c_{X\ot Y, Z}=
\gbeg{3}{4}
\got{1}{X}\got{1}{Y}\got{1}{Z}\gnl
\gcl{1}\gbr\gnl
\gbr\gcl{1}\gnl
\gob{1}{Z}\gob{1}{X}\gob{1}{Y}
\gend
\hspace{2mm},
\end{equation}
for all objects $X, Y, Z\in \Cc$, and 
\begin{equation}\eqlabel{qybe}
\gbeg{3}{5}
\got{1}{X}\got{1}{Y}\got{1}{Z}\gnl
\gbr\gcl{1}\gnl
\gcl{1}\gbr\gnl
\gbr\gcl{1}\gnl
\gob{1}{Z}\gob{1}{Y}\gob{1}{X}
\gend =
\gbeg{3}{5}
\got{1}{X}\got{1}{Y}\got{1}{Z}\gnl
\gcl{1}\gbr\gnl
\gbr\gcl{1}\gnl
\gcl{1}\gbr\gnl
\gob{1}{Z}\gob{1}{Y}\gob{1}{X}
\gend
\hspace{1mm},
\end{equation}
the categorical version of the Yang-Baxter equation. Furthermore, $c$ 
is natural in the sense that 
\[
\gbeg{2}{4}
\got{1}{M}\got{1}{N}\gnl
\gbr\gnl
\gmp{g}\gmp{f}\gnl
\gob{1}{V}\gob{1}{U}
\gend =
\gbeg{2}{4}
\got{1}{M}\got{1}{N}\gnl
\gmp{f}\gmp{g}\gnl
\gbr\gnl
\gob{1}{V}\gob{1}{U}
\gend 
\hspace{2mm},    
\]
for all morphisms $f:M\ra U$ and $g: N\ra V$ in ${\cal C}$. 
In particular, if we have a morphism $X\ot Y\to Z$ in ${\cal C}$,
which is often denoted by 
$
\gbeg{3}{3} 
\got{1}{X}\gvac{1}\got{1}{Y}\gnl
\gwmu{3}\gnl
\gvac{1}\gob{1}{Z}
\gend
$, then we have 
\begin{equation}\eqlabel{nat1cup}
\gbeg{4}{5}
\got{1}{T}\got{1}{X}\gvac{1}\got{1}{Y}\gnl
\gcl{1}\gwmu{3}\gnl\
\gcl{1}\gcn{1}{1}{3}{1}\gnl
\gbr\gnl
\gob{1}{Z}\gob{1}{T}
\gend
=
\gbeg{4}{5}
\gvac{1}\got{1}{T}\got{1}{X}\got{1}{Y}\gnl
\gvac{1}\gbr\gcl{1}\gnl
\gvac{1}\gcn{1}{1}{1}{-1}\gbr\gnl
\gwmu{3}\gcl{1}\gnl
\gvac{1}\gob{1}{Z}\gvac{1}\gob{1}{T}
\gend
\hspace{2mm}\mbox{and}\hspace{2mm}
\gbeg{4}{5}
\got{1}{X}\gvac{1}\got{1}{Y}\got{1}{T}\gnl
\gwmu{3}\gcl{1}\gnl
\gvac{1}\gcn{1}{1}{1}{3}\gvac{1}\gcl{1}\gnl
\gvac{2}\gbr\gnl
\gvac{2}\gob{1}{T}\gob{1}{Z}\gnl
\gend 
=
\gbeg{4}{5}
\got{1}{X}\got{1}{Y}\got{1}{T}\gnl
\gcl{1}\gbr\gnl
\gbr\gcn{1}{1}{1}{3}\gnl
\gcl{1}\gwmu{3}\gnl
\gob{1}{T}\gvac{1}\gob{1}{Z}
\gend
\hspace{2mm},
\end{equation}
for all $T\in {\cal C}$. Similarly, if 
$
\gbeg{3}{3}
\gvac{1}\got{1}{X}\gnl
\gwcm{3}\gnl
\gob{1}{Y}\gvac{1}\gob{1}{Z}\gnl
\gend
$ 
is a morphism from $X$ to $Y\ot Z$ in $\Cc$ then the naturality 
of $c$ implies that
\begin{equation}\eqlabel{nat2cup}
\gbeg{4}{5}
\got{1}{X}\got{1}{T}\gnl
\gbr\gnl
\gcl{1}\gcn{1}{1}{1}{3}\gnl
\gcl{1}\gwcm{3}\gnl
\gob{1}{T}\gob{1}{Y}\gvac{1}\gob{1}{Z}
\gend
=
\gbeg{4}{5}
\gvac{1}\got{1}{X}\gvac{1}\got{1}{T}\gnl
\gwcm{3}\gcl{1}\gnl
\gcn{1}{1}{1}{3}\gvac{1}\gbr\gnl
\gvac{1}\gbr\gcl{1}\gnl
\gvac{1}\gob{1}{T}\gob{1}{Y}\gob{1}{Z}
\gend
\hspace{2mm}\mbox{and}\hspace{2mm}
\gbeg{4}{4}
\gvac{1}\got{1}{T}\got{1}{X}\gnl
\gvac{1}\gbr\gnl
\gwcm{3}\gcn{1}{1}{-1}{1}\gnl
\gob{1}{Y}\gvac{1}\gob{1}{Z}\gob{1}{T}
\gend 
=
\gbeg{4}{5}
\got{1}{T}\gvac{1}\got{1}{X}\gnl
\gcl{1}\gwcm{3}\gnl
\gbr\gvac{1}\gcn{1}{1}{1}{-1}\gnl
\gcl{1}\gbr\gnl
\gob{1}{Y}\gob{1}{Z}\gob{1}{T}
\gend
\hspace{2mm},
\end{equation}
for all $T\in {\cal C}$. For $X\in \Cc$, we identify $\un{1}\ot X\cong X\cong X\ot\un{1}$ using the 
left and right unit constraints. By \cite[Prop. XIII.1.2]{kas} 
we can also identify $c_{\un{1}, X}$ and $c_{X, \un{1}}$ with the identity morphism 
of $X$ in $\Cc$, which will be denoted from now on by 
$
\Id_X=
\gbeg{1}{3}
\got{1}{X}\gnl
\gcl{1}\gnl
\gob{1}{X}
\gend
$. 
In addition, all the results will be proved for strict monoidal categories, 
these are monoidal categories for which all the associativity, left and right unit 
constraints are identity morphisms. The results remain valid for an arbitrary
monoidal category, since every monoidal category is equivalent to a strict one,
see \cite{kas}.

\subsection{Braided bialgebras}\selabel{1.2}
An algebra in a monoidal category $\Cc$ is an object $A$ of 
$\Cc$ endowed with a multiplication $\un{m}_A:\ A\ot A\ra A$ and a unit morphism 
$\un{\eta}_A:\ \un{1}\ra A$ which are associative and unital up to the
associativity and unit constraints. 
The multiplication and the unit of $A$ will be denoted by 
$
\gbeg{2}{3} 
\got{1}{A}\got{1}{A}\gnl
\gmu\gnl
\gob{2}{A}
\gend
$ 
and 
$\gbeg{1}{3}
\got{1}{\un{1}}\gnl
\gu{1}\gnl
\gob{1}{A}
\gend
$.
A coalgebra $B$ in $\Cc$ is an algebra in the opposite category. The comultiplication
$\un{\Delta}_B:\ B\ra B\ot B$ and counit $\un{\va}_B: B\ra \un{1}$ will respectively be
denoted by 
$
\gbeg{2}{3}
\got{2}{B}\gnl
\gcmu\gnl
\gob{1}{B}\gob{1}{B}
\gend
$ 
and 
$
\gbeg{1}{3}
\got{1}{B}\gnl
\gcu{1}\gnl
\gob{1}{\un{1}}
\gend
$.\\
A bialgebra $(B, \un{m}_B, \un{\eta }_B, \un {\Delta }_B, \un {\va }_B)$ in a prebraided 
monoidal category $\Cc$ is a 5-tuple $(B, \un{m}_B, \un{\eta }_B, \un {\Delta }_B, \un {\va }_B)$,
with $(B, \un{m}_B, \un{\eta}_B)$ an algebra and $(B, \un{\Delta}_B, \un{\va}_B)$
a coalgebra in $\Cc$ such that
$\un{\Delta}_B:\ B\to B\ot B$ and $\un{\va}_B: B\ra \un{1}$ are algebra morphisms. 
Here $B\ot B$ has the tensor product algebra structure and $\un{1}$ is viewed as an algebra 
in $\Cc$ through the left or right unit constraint. Explicitly, apart from $\un{\va}_B\un{\eta}_B=\Id_{\un{1}}$, 
in diagrammatic notations the axioms for a bialgebra $B$ in $\Cc$ read as follows, 
\begin{eqnarray}
&&
\gbeg{3}{5}
\got{1}{B}\got{1}{B}\got{1}{B}\gnl
\gmu\gcl{2}\gnl
\gvac{1}\gcn{1}{1}{0}{1}\gnl
\gvac{1}\gmu\gnl
\gob{4}{B}
\gend =
\gbeg{3}{5}
\got{1}{B}\got{1}{B}\got{1}{B}\gnl
\gcl{2}\gmu\gnl
\gvac{1}\gcn{1}{1}{2}{1}\gnl
\gmu\gnl
\gob{2}{B}
\gend
\hspace{2mm},\hspace{2mm}
\gbeg{2}{4}
\got{1}{B}\gnl
\gcl{1}\gu{1}\gnl
\gmu\gnl
\gob{2}{B}
\gend =
\gbeg{1}{3}
\got{1}{B}\gnl
\gcl{1}\gnl
\gob{1}{B}
\gend =
\gbeg{2}{4}
\got{3}{B}\gnl
\gu{1}\gcl{1}\gnl
\gmu\gnl
\gob{2}{B}
\gend
\hspace{2mm},\hspace{2mm}
\gbeg{3}{5}
\got{4}{B}\gnl
\gvac{1}\gcmu\gnl
\gvac{1}\gcn{1}{1}{1}{0}\gcl{1}\gnl
\gcmu\gcl{1}\gnl
\gob{1}{B}\gob{1}B\gob{1}{B}
\gend =
\gbeg{3}{5}
\got{2}{B}\gnl
\gcmu\gnl
\gcl{2}\gcn{1}{1}{1}{2}\gnl
\gvac{1}\gcmu\gnl
\gob{1}B\gob{1}{B}\gob{1}{B}
\gend 
\hspace{2mm},
\nonumber\\
&&\label{braidedbialgebra}\\
&&
\gbeg{2}{4}
\got{2}{B}\gnl
\gcmu\gnl
\gcl{1}\gcu{1}\gnl
\gob{1}{B}
\gend =
\gbeg{1}{3}
\got{1}{B}\gnl
\gcl{1}\gnl
\gob{1}{B}
\gend =
\gbeg{2}{4}
\got{2}{B}\gnl
\gcmu\gnl
\gcu{1}\gcl{1}\gnl
\gob{3}{B}
\gend
\hspace{1mm},\hspace{1mm}
\gbeg{2}{5}
\got{1}{B}\got{1}{B}\gnl
\gmu\gnl
\gcn{1}{1}{2}{1}\gnl
\gcu{1}\gnl
\gob{2}{\un{1}}
\gend
=
\gbeg{2}{3}
\got{1}{B}\got{1}{B}\gnl
\gcu{1}\gcu{1}\gnl
\gob{2}{\un{1}}
\gend
\hspace{1mm},\hspace{1mm}
\gbeg{2}{5}
\got{2}{\un{1}}\gnl
\gu{1}\gnl
\gcn{1}{1}{1}{2}\gnl
\gcmu\gnl
\gob{1}{B}\gob{1}{B}
\gend
=
\gbeg{2}{3}
\got{2}{\un{1}}\gnl
\gu{1}\gu{1}\gnl
\gob{1}{B}\gob{1}{B}
\gend
\hspace{1mm},\hspace{1mm}
\gbeg{2}{4}
\got{1}{B}\got{1}{B}\gnl
\gmu\gnl
\gcmu\gnl
\gob{1}{B}\gob{1}{B}
\gend = 
\gbeg{4}{5}
\got{2}{B}\got{2}{B}\gnl
\gcmu\gcmu\gnl
\gcl{1}\gbr\gcl{1}\gnl
\gmu\gmu\gnl
\gob{2}{B}\gob{2}{B}
\gend
\mbox{\hspace{2mm}.}
\nonumber
\end{eqnarray}
If $B$ is a bialgebra in $\Cc$ then ${}_B\Cc$ (resp. ${}^B\Cc$), the category of left 
$B$-modules (resp. left $B$-comodules) in $\Cc$ is a monoidal category. If $X, Y$ 
are objects in ${}_B\Cc$ (resp. ${}^B\Cc$) then $X\ot Y$ is a left $B$-module (resp. 
left $B$-comodule) via the action (resp. coaction)
\begin{equation}\eqlabel{monstrforrepcorep}
\gbeg{4}{5}
\got{2}{B}\got{1}{X}\got{1}{Y}\gnl
\gcmu\gcl{1}\gcl{1}\gnl
\gcl{1}\gbr\gcl{1}\gnl
\glm\glm\gnl
\gvac{1}\gob{1}{X}\gvac{1}\gob{1}{Y}
\gend
\hspace{3mm}
\left(\mbox{resp.}
\hspace{3mm}
\gbeg{4}{5}
\gvac{1}\got{1}{X}\gvac{1}\got{1}{Y}\gnl
\glcm\glcm\gnl
\gcl{1}\gbr\gcl{1}\gnl
\gmu\gcl{1}\gcl{1}\gnl
\gob{2}{B}\gob{1}{X}\gob{1}{Y}
\gend
\right)
\hspace{2mm},
\end{equation}
where 
$
\gbeg{2}{3}
\got{1}{B}\got{1}{X}\gnl
\glm\gnl
\gvac{1}\gob{1}{X}
\gend
$ 
is our diagrammatic notation for the left action of $B$ on $X$, while 
$
\gbeg{2}{3}
\gvac{1}\got{1}{X}\gnl
\glcm\gnl
\gob{1}{B}\gob{1}{X}
\gend
$ 
is the notation used for the left $B$-coaction on $X$, etc.\\
If $B$ is a braided bialgebra, then we can consider algebras and coalgebras in
${}^B\Cc$ and ${}_B{\cal C}$. 
A (co)algebra in ${}^B\Cc$ is called a left $B$-comodule (co)algebra, and a
(co)algebra in ${}_B\Cc$ is called a left $B$-module (co)algebra.
More precisely, with notation as above, a 
left $B$-comodule (co)algebra in $\Cc$ is a left $B$-comodule $A$ which is at the same time a (co)algebra in $\Cc$ such that 
the (co)multiplication and the (co)unit are morphisms in ${}^B\Cc$, that is,
\begin{equation}\eqlabel{lcomcoalg}
\gbeg{2}{4}
\gvac{1}\got{1}{A}\gnl
\glcm\gnl
\gcl{1}\gcu{1}\gnl
\gob{1}{B}
\gend
=
\gbeg{2}{3}
\gvac{1}\got{1}{A}\gnl
\gu{1}\gcu{1}\gnl
\gob{1}{B}
\gend
\hspace{2mm},\hspace{2mm}
\gbeg{4}{7}
\gvac{1}\got{2}{A}\gnl
\gvac{1}\gcmu\gnl
\glcm\gcn{1}{1}{1}{3}\gnl
\gcl{1}\gcl{1}\glcm\gnl
\gcl{1}\gbr\gcl{1}\gnl
\gmu\gcl{1}\gcl{1}\gnl
\gob{2}{B}\gob{1}{A}\gob{1}{A}
\gend
=
\gbeg{3}{5}
\gvac{1}\got{1}{A}\gnl
\glcm\gnl
\gcl{1}\gcn{1}{1}{1}{2}\gnl
\gcl{1}\gcmu\gnl
\gob{1}{B}\gob{1}{A}\gob{1}{A}
\gend
\end{equation}
in the comodule coalgebra case, and
\begin{equation}\eqlabel{lcomalg}
\gbeg{2}{4}
\got{2}{\un{1}}\gnl
\gvac{1}\gu{1}\gnl
\glcm\gnl
\gob{1}{B}\gob{1}{A}
\gend
=
\gbeg{2}{3}
\got{2}{\un{1}}\gnl
\gu{1}\gu{1}\gnl
\gob{1}{B}\gob{1}{A}
\gend
\hspace{2mm},\hspace{2mm}
\gbeg{3}{5}
\gvac{1}\got{1}{A}\got{1}{A}\gnl
\gvac{1}\gmu\gnl
\gvac{1}\gcn{1}{1}{2}{1}\gnl
\glcm\gnl
\gob{1}{B}\gob{1}{A}
\gend
=
\gbeg{4}{5}
\gvac{1}\got{1}{A}\gvac{1}\got{1}{A}\gnl
\glcm\glcm\gnl
\gcl{1}\gbr\gcl{1}\gnl
\gmu\gmu\gnl
\gob{2}{B}\gob{2}{A}
\gend
\end{equation}
in the comodule algebra case.\\
In a similar way, a (co)algebra in ${}_B\Cc$ is called a left $B$-module (co)algebra,
that is a left $B$-module with a (co)algebra structure in $\Cc$ such that the 
(co)multiplication and (co)unit are morphisms in ${}_B\Cc$; this condition can be
expressed in a diagrammatic way as follows:
\begin{equation}\eqlabel{lmodalg}
\gbeg{2}{4}
\got{1}{B}\gnl
\gcl{1}\gu{1}\gnl
\glm\gnl
\gvac{1}\gob{1}{A}
\gend
=
\gbeg{2}{3}
\got{1}{B}\gnl
\gcu{1}\gu{1}\gnl
\gvac{1}\gob{1}{A}
\gend
\hspace{2mm},\hspace{2mm}
\gbeg{3}{5}
\got{1}{B}\got{1}{A}\got{1}{A}\gnl
\gcl{1}\gmu\gnl
\gcl{1}\gcn{1}{1}{2}{1}\gnl
\glm\gnl
\gvac{1}\gob{1}{A}
\gend
=
\gbeg{4}{7}
\got{2}{B}\got{1}{A}\got{1}{A}\gnl
\gcmu\gcl{1}\gcl{1}\gnl
\gcl{1}\gbr\gcl{1}\gnl
\glm\glm\gnl
\gvac{1}\gcl{1}\gvac{1}\gcn{1}{1}{1}{-1}\gnl
\gvac{1}\gmu\gnl
\gvac{1}\gob{2}{A}
\gend
\hspace{2mm},
\end{equation}
in the algebra case, and
\begin{equation}\eqlabel{lmodcoalg}
\gbeg{2}{4}
\got{1}{B}\got{1}{A}\gnl
\glm\gnl
\gvac{1}\gcu{1}\gnl
\gob{2}{\un{1}}
\gend
=
\gbeg{2}{3}
\got{1}{B}\got{1}{A}\gnl
\gcu{1}\gcu{1}\gnl
\gob{2}{\un{1}}
\gend
\hspace{2mm},\hspace{2mm}
\gbeg{3}{5}
\got{1}{B}\got{1}{A}\gnl
\glm\gnl
\gvac{1}\gcn{1}{1}{1}{0}\gnl
\gcmu\gnl
\gob{1}{A}\gob{1}{A}
\gend
=
\gbeg{4}{5}
\got{2}{B}\got{2}{A}\gnl
\gcmu\gcmu\gnl
\gcl{1}\gbr\gcl{1}\gnl
\glm\glm\gnl
\gvac{1}\gob{1}{A}\gvac{1}\gob{1}{A}
\gend
\hspace{2mm}.
\end{equation}
in the coalgebra case.

\subsection{Yetter-Drinfeld modules}\selabel{1.3}
The category of left Yetter-Drinfeld modules ${}_B^B{\cal YD}$  over a bialgebra $B$ in a braided
monoidal category $\Cc$ was introduced in \cite{besp}.
It is a prebraided monoidal category, and 
it can be identified with a full subcategory of the left weak center of the monoidal 
category ${}_B\Cc$ (see \cite[Prop. 3.6.1]{besp}). We now give an explicit description.\\
A left Yetter-Drinfeld module is an object $X\in \Cc$ with a left $B$-action and a left
$B$-coaction satisfying the compatibility relation
\begin{equation}\eqlabel{lYDMcond} 
\gbeg{3}{9}
\got{2}{B}\got{1}{X}\gnl
\gcmu\gcl{1}\gnl
\gcl{1}\gbr\gnl
\gev\gvac{-1}\gcl{2}\gcl{1}\gnl
\gvac{1}\gcl{2}\gcl{1}\gnl
\gdb\gcl{1}\gnl
\gcl{1}\gbr\gnl
\gmu\gcl{1}\gnl
\gob{2}{B}\gob{1}{X}
\gend
=
\gbeg{4}{7}
\got{2}{B}\gvac{1}\got{1}{X}\gnl
\gcmu\gvac{1}\gcl{2}\gnl
\gcl{1}\gcl{1}\gdb\gnl
\gcl{1}\gbr\gcl{1}\gnl
\gmu\gev\gvac{-1}\gcl{2}\gnl
\gcn{1}{1}{2}{2}\gnl
\gob{2}{B}\gvac{1}\gob{1}{X}
\gend
\hspace{2mm},
\end{equation}
where, from now on and in order to avoid confusion, we denote a $B$-action on a generic $X$, different from 
$A$ and $B$, by 
$
\gbeg{2}{4}
\got{1}{B}\got{1}{X}\gnl
\gev\gvac{-1}\gcl{2}\gnl
\gvac{1}\gob{1}{X}
\gend
$;
similarly, for an object $X$ different from $A, B$ we denote a left coaction of $B$ on $X$ by 
$
\gbeg{2}{4}
\gvac{1}\got{1}{X}\gnl
\gvac{1}\gcl{2}\gnl
\gdb\gnl
\gob{1}{B}\gob{1}{X}
\gend
$. \\
Morphisms in ${}^B_B{\cal YD}$ are morphisms in $\Cc$ that are left $B$-linear and left
$B$-colinear.\\
The tensor and prebraiding on ${}^B_B{\cal YD}$ are inherited from the tensor and prebraiding 
on the left weak center of ${}_B\Cc$. Namely, the left $B$-action and left $B$-coaction
on the tensor product $X\ot Y$ of $X,Y\in {}^B_B{\cal YD}$ is given by 
\equref{monstrforrepcorep}, and the prebraiding $\un{c}$ is defined by
\[
\un{c}_{X, Y}=
\gbeg{3}{7}
\gvac{1}\got{1}{X}\got{1}{Y}\gnl
\gvac{1}\gcl{2}\gcl{1}\gnl
\gdb\gcl{1}\gnl
\gcl{1}\gbr\gnl
\gev\gvac{-1}\gcl{2}\gcl{2}\gnl
\gvac{1}\gob{1}{Y}\gob{1}{X}
\gend
\hspace{2mm}.
\] 
Consequently, a bialgebra in ${}_B^B{\cal YD}$ is an object $A$ in ${}_B^B{\cal YD}$ that admits an algebra and a coalgebra 
structure in $\Cc$ satisfying conditions (\ref{eq:lcomcoalg}-\ref{eq:lmodcoalg}) and
\begin{equation}\eqlabel{bialgYD}
\un{\va}_A\un{\eta}_A=\Id_{\un{1}}
\hspace{1mm},\hspace{1mm}
\gbeg{2}{5}
\got{1}{A}\got{1}{A}\gnl
\gmu\gnl
\gcn{1}{1}{2}{1}\gnl
\gcu{1}\gnl
\gob{2}{\un{1}}
\gend
=
\gbeg{2}{3}
\got{1}{A}\got{1}{A}\gnl
\gcu{1}\gcu{1}\gnl
\gob{2}{\un{1}}
\gend
\hspace{1mm},\hspace{1mm}
\gbeg{2}{5}
\got{2}{\un{1}}\gnl
\gu{1}\gnl
\gcn{1}{1}{1}{2}\gnl
\gcmu\gnl
\gob{1}{A}\gob{1}{A}
\gend
=
\gbeg{2}{3}
\got{2}{\un{1}}\gnl
\gu{1}\gu{1}\gnl
\gob{1}{A}\gob{1}{A}
\gend
\hspace{1mm},\hspace{1mm}
\gbeg{2}{4}
\got{1}{A}\got{1}{A}\gnl
\gmu\gnl
\gcmu\gnl
\gob{1}{A}\gob{1}{A}
\gend
=
\gbeg{5}{9}
\got{2}{A}\gvac{1}\got{2}{A}\gnl
\gcmu\gvac{1}\gcmu\gnl
\gcl{1}\gcn{1}{1}{1}{3}\gvac{1}\gcl{1}\gcl{1}\gnl
\gcl{1}\glcm\gcl{1}\gcl{1}\gnl
\gcl{1}\gcl{1}\gbr\gcl{1}\gnl
\gcl{1}\glm\gmu\gnl
\gcl{1}\gcn{1}{1}{3}{1}\gvac{1}\gcn{1}{2}{2}{2}\gnl
\gmu\gnl
\gob{2}{A}\gvac{1}\gob{2}{A}
\gend
\hspace{1mm},
\end{equation}
where we used the notation from \seref{1.2}. Note that the above equations express the fact that the counit and the comultiplication of 
$A$ are algebra morphisms in ${}_B^B{\cal YD}$. 

\section{Main result}\selabel{2}
\setcounter{equation}{0}

From now on, we assume that $B$ is a bialgebra in a braided monoidal category $\Cc$,
and that $A$ is a left 
$B$-comodule algebra in $\Cc$ with left $B$-coaction 
$
\gbeg{2}{3}
\gvac{1}\got{1}{A}\gnl
\glcm\gnl
\gob{1}{B}\gob{1}{A}
\gend
$. 
${}^B\Cc_A$ will be the notation for the category of right-left relative 
$(B, A)$-Hopf modules in $\Cc$. These are objects $X\in \Cc$ with a left $B$-coaction and
a right $A$-action such that 
\begin{equation}\eqlabel{rDHcond}
\gbeg{3}{5}
\gvac{1}\got{1}{X}\got{1}{A}\gnl
\gvac{1}\gev\gvac{-2}\gcl{2}\gnl
\gvac{1}\gcl{1}\gnl
\gdb\gvac{-1}\gcl{1}\gnl
\gob{1}{B}\gob{1}{X}
\gend
=
\gbeg{4}{7}
\gvac{1}\got{1}{X}\gvac{1}\got{1}{A}\gnl
\gvac{1}\gcl{2}\glcm\gnl
\gdb\gcl{1}\gcl{1}\gnl
\gcl{1}\gbr\gcl{1}\gnl
\gmu\gcl{2}\gcl{1}\gnl
\gcn{1}{1}{2}{2}\gvac{1}\gev\gnl
\gob{2}{B}\gob{1}{X}
\gend
\hspace{2mm}.
\end{equation}
Morphisms in ${}^B\Cc_A$ are morphisms in $\Cc$ that are left $B$-colinear 
and right $A$-linear.\\
Observe that $A\in {}^B\Cc_A$, with right $A$-action given by multiplication and
left $B$-coaction via the $B$-comodule algebra structure.\\
We call $(B,A)$ an input monoidal Doi-Hopf datum if $B$ is a bialgebra, $A$ is a
left $B$-comodule algebra and a coalgebra, and we have a morphism $B\ot A\to A$.
We do not assume that this makes $A$ into bialgebra or a left $B$-module coalgebra. Nevertheless,
we still use the diagrammatic notation
$
\gbeg{2}{3}
\got{1}{B}\got{1}{A}\gnl
\glm\gnl
\gvac{1}\gob{1}{A}
\gend
$
for the (not necessarily associative or unital) $B$-action $B\ot A\to A$.\\
Now take two relative Hopf modules $X$ and $Y$. We know that $X\ot Y$ is a left
$B$-comodule using \equref{monstrforrepcorep}. Assume that $(B,A)$ is an input monoidal Doi-Hopf datum.
\begin{itemize}
\item We can define a right $A$-action on $X\ot Y$ using the diagram
\begin{equation}\eqlabel{rAactDHM}
\gbeg{5}{9}  
\got{1}{X}\gvac{1}\got{1}{Y}\got{2}{A}\gnl
\gcl{1}\gvac{1}\gcl{2}\gcmu\gnl
\gcl{1}\gdb\gcl{1}\gcl{1}\gnl
\gcl{1}\gcl{1}\gbr\gcl{1}\gnl
\gcl{1}\glm\gcl{4}\gvac{-1}\gev\gnl
\gcl{1}\gcn{1}{1}{3}{1}\gnl
\gcl{1}\gvac{-1}\gev\gnl
\gcl{1}\gnl
\gob{1}{X}\gvac{2}\gob{1}{Y}
\gend
\end{equation}
\item $\un{1}$ has a left $B$-coaction defined by the unit of $B$ and a right $A$-action defined
by the counit of $A$.
\item For $X,Y,Z\in {}^B\Cc_A$, we can write down the associativity constraints $a_{X,Y,Z}$
and the unit constraints $l_X$ and $r_X$ in $\Cc$.
\end{itemize}
This provides part of the ingredients that are needed to define a monoidal structure on the
category of relative Hopf Modules ${}^B\Cc_A$. We will say that 
$({}^B\Cc_A, \ot, \un{1}, a, l, r)$ is the input monoidal structure on ${}^B\Cc_A$
defined by the input monoidal Hopf module datum $(B, A)$. The main result of this
note is the following.

\begin{theorem}\thlabel{2.1}
The input monoidal structure on ${}^B\Cc_A$
defined by the input monoidal Hopf module datum $(B, A)$ is a monoidal structure if and
only if $A$ is a bialgebra in the prebraided monoidal category ${}^B_B{\cal YD}$.
\end{theorem}

Before we present the proof of \thref{2.1}, we need some Lemmas.

\begin{lemma}\lelabel{unitcounitcomp}
Let $({}^B\Cc_A, \ot, \un{1}, a, l, r)$ be the input monoidal structure on ${}^B\Cc_A$ associated to an input monoidal Hopf module 
datum $(B, A)$. The following statements are equivalent.\\
1) $\un{1}\in {}^B\Cc_A$, $l_X$ and $r_X$ are morphisms in $ {}^B\Cc_A$,
and the tensor product $X\ot Y$ satisfies the unit condition, for all $X,Y\in {}^B\Cc_A$.\\
2) we have the following compatibility relations between
the unit and counit morphisms of $A$ and $B$,
\begin{eqnarray}
&&\gbeg{2}{5}
\got{1}{A}\got{1}{A}\gnl
\gmu\gnl
\gcn{1}{1}{2}{1}\gnl
\gcu{1}\gnl
\gob{2}{\un{1}}
\gend
=
\gbeg{2}{3}
\got{1}{A}\got{1}{A}\gnl
\gcu{1}\gcu{1}\gnl
\gob{2}{\un{1}}
\gend
\hspace{1mm},\hspace{1mm}
\un{\va}_A\un{\eta}_A=\Id_{\un{1}}
\hspace{1mm},\hspace{1mm}
\gbeg{2}{4}
\gvac{1}\got{1}{A}\gnl
\glcm\gnl
\gcl{1}\gcu{1}\gnl
\gob{1}{B}
\gend
=
\gbeg{2}{3}
\gvac{1}\got{1}{A}\gnl
\gu{1}\gcu{1}\gnl
\gob{1}{B}
\gend
\hspace{1mm},\hspace{1mm}
\gbeg{2}{4}
\got{1}{B}\got{1}{A}\gnl
\glm\gnl
\gvac{1}\gcu{1}\gnl
\gob{2}{\un{1}}
\gend
=
\gbeg{2}{3}
\got{1}{B}\got{1}{A}\gnl
\gcu{1}\gcu{1}\gnl
\gob{2}{\un{1}}
\gend
\hspace{1mm},\hspace{1mm}\nonumber\\
&&\eqlabel{unitcounitrel}\\
&&
\gbeg{2}{4}
\got{1}{B}\got{1}{A}\gnl
\gu{1}\gcl{1}\gnl
\glm\gnl
\gvac{1}\gob{1}{A}
\gend
=
\gbeg{1}{3}
\got{1}{A}\gnl
\gcl{1}\gnl
\gob{1}{A}
\gend
\hspace{1mm},\hspace{1mm}
\gbeg{2}{5}
\got{2}{\un{1}}\gnl
\gu{1}\gnl
\gcn{1}{1}{1}{2}\gnl
\gcmu\gnl
\gob{1}{A}\gob{1}{A}\gnl
\gend
=
\gbeg{2}{3}
\got{2}{\un{1}}\gnl
\gu{1}\gu{1}\gnl
\gob{1}{A}\gob{1}{A}
\gend
\hspace{1mm},\hspace{1mm}
\gbeg{2}{4}
\got{1}{B}\gnl
\gcl{1}\gu{1}\gnl
\glm\gnl
\gvac{1}\gob{1}{A}
\gend
=
\gbeg{2}{4}
\got{1}{B}\gnl
\gcu{1}\gnl
\gvac{1}\gu{1}\gnl
\gvac{1}\gob{1}{A}
\gend
\hspace{1mm}.
\nonumber
\end{eqnarray}
In particular, if the input  monoidal structure is monoidal, then the compatility relations
\equref{unitcounitrel} hold.
\end{lemma}

\begin{proof}
We examine first when $\un{1}$ is an object of ${}^B\Cc_A$. Since $\un{\Delta}_B$ and $\un{\va}_B$ respect the 
unit of $B$ it follows that $\un{1}$ is always a left $B$-comodule in $\Cc$ via the unit morphism of $B$. Now, it 
can be easily checked that $\un{1}$ is a right $A$-module in $\Cc$ via the counit of $A$ if and only if $\un{\va}_A$ is an 
algebra morphism in $\Cc$, and that the module-comodule compatibilty relation holds in this case if and only if 
$
\gbeg{2}{4}
\gvac{1}\got{1}{A}\gnl
\glcm\gnl
\gcl{1}\gcu{1}\gnl
\gob{1}{B}
\gend
=
\gbeg{2}{3}
\gvac{1}\got{1}{A}\gnl
\gu{1}\gcu{1}\gnl
\gob{1}{B}
\gend
$. 
Thus we have shown that $\un{1}\in {}^B\Cc_A$ if and only if
the first three equalities in \equref{unitcounitrel} hold.\\
We know at this moment that $\un{\va}_A$ is an algebra morphism in $\Cc$,
hence $\un{\eta}_B\circ \un{\va}_A:\ A\to B$ is also an algebra morphism, and
$B$ can be viewed as a left $A$-module via restriction of scalars.
$B$ is also a left $B$-comodule in $\Cc$ via its comultiplication and it follows from the third equality in 
\equref{unitcounitrel} that $B\in {}^B\Cc_A$. We will call this relative Hopf module structure
on $B$ trivial, and denote $B$ with this structure by $B_{\rm tr}$.\\
$l_X$ is always left $B$-colinear; this follows from a simple inspection, and is due to the fact
that the category ${}^B\Cc$ is monoidal. $l_X$ is right $A$-linear if and only if
\[
\gbeg{4}{7}
\gvac{1}\got{1}{X}\got{2}{A}\gnl
\gvac{1}\gcl{2}\gcmu\gnl
\gdb\gcl{1}\gcl{1}\gnl
\gcl{1}\gbr\gcl{1}\gnl
\glm\gcl{2}\gvac{-1}\gev\gnl
\gvac{1}\gcu{1}\gnl
\gvac{2}\gob{1}{X}
\gend
=
\gbeg{2}{4}
\got{1}{X}\got{1}{A}\gnl
\gcl{2}\gvac{-1}\gev\gnl
\gob{1}{X}
\gend
\hspace{1mm}.
\]
We conclude that the left unit constraint morphisms  $l_X$ of relative Hopf modules $X$
are in ${}^B\Cc_A$ if and only if the fourth relation in \equref{unitcounitrel} holds.
To see the direct implication, take $X=B_{\rm tr}$ and then apply $\un{\va}_B$ to the lower $B$.\\
In a similar way, $r_X$ is always left $B$-colinear, and is right $A$-linear if and only if
\[
\gbeg{3}{7}
\got{1}{X}\gvac{1}\got{1}{A}\gnl
\gcl{1}\gu{1}\gcl{1}\gnl
\gcl{1}\glm\gnl
\gcl{1}\gcn{1}{1}{3}{1}\gnl
\gcl{2}\gvac{-1}\gev\gnl
\gob{1}{X}
\gend
=
\gbeg{2}{4}
\got{1}{X}\got{1}{A}\gnl
\gcl{2}\gvac{-1}\gev\gnl
\gob{1}{X}
\gend
\hspace{1mm}.
\]
Then it follows that the right unit constraint morphisms of all relative Hopf modules $X$
are in ${}^B\Cc_A$ if and only if the fifth relation in \equref{unitcounitrel} holds.
For the direct implication, take $X=A$ in the above equality, and then compose it to the right
by $\un{\eta}_A\ot {\rm id}_A$.\\
We are left to show the two final equalities in \equref{unitcounitrel}. We will see that they
follow from the unit condition on the tensor product of $X,Y\in {}^B\Cc_A$. More precisely,
the right $A$-module structure on $X\ot Y$ respects the unit of $A$ if and only if
\begin{equation}\eqlabel{unitcondtensprod}
\gbeg{5}{10}
\got{1}{X}\gvac{1}\got{1}{Y}\gnl
\gcl{1}\gvac{1}\gcl{2}\gu{1}\gnl
\gcl{1}\gdb\gcn{1}{1}{1}{2}\gnl
\gcl{1}\gcl{1}\gcl{1}\gcmu\gnl
\gcl{1}\gcl{1}\gbr\gcl{1}\gnl
\gcl{1}\glm\gcl{4}\gvac{-1}\gev\gnl
\gcl{1}\gcn{1}{1}{3}{1}\gnl
\gev\gvac{-2}\gcl{2}\gnl
\gob{1}{X}\gvac{2}\gob{1}{Y}
\gend
=
\gbeg{2}{3}
\got{1}{X}\got{1}{Y}\gnl
\gcl{1}\gcl{1}\gnl
\gob{1}{X}\gob{1}{Y}
\gend
\hspace*{5mm}\hbox{if and only if}\hspace*{5mm}
\gbeg{4}{8}
\gvac{1}\got{1}{Y}\gnl
\gvac{1}\gcl{3}\gu{1}\gnl
\gdb\gcn{1}{1}{1}{2}\gnl
\gcl{1}\gvac{1}\gcmu\gnl
\gcl{1}\gbr\gcl{1}\gnl
\glm\gcl{2}\gvac{-1}\gev\gnl
\gvac{1}\gcl{1}\gnl
\gvac{1}\gob{1}{A}\gob{1}{Y}
\gend
=
\gbeg{2}{3}
\gvac{1}\got{1}{Y}\gnl
\gu{1}\gcl{1}\gnl
\gob{1}{A}\gob{1}{Y}
\gend
\hspace{1mm}.
\end{equation}
For the direct implication, take $X=A$, and compose the equality to the right by 
$\un{\eta}_A\ot {\rm id}_Y$. Consequently, the unit condition on $X\ot Y$ is equivalent to
\equref{unitcondtensprod}. Now take $Y=A$ in the second equality of \equref{unitcondtensprod},
and compose it to the right by $\un{\eta}_A$. With the help of 
$
\gbeg{2}{4}
\got{2}{\un{1}}\gnl
\gvac{1}\gu{1}\gnl
\glcm\gnl
\gob{1}{B}\gob{1}{A}
\gend
 =
\gbeg{2}{3}
\got{2}{\un{1}}\gnl
\gu{1}\gu{1}\gnl
\gob{1}{B}\gob{1}{A}
\gend
$ 
and 
$
\gbeg{2}{4}
\gvac{1}\got{1}{A}\gnl
\gu{1}\gcl{1}\gnl
\glm\gnl
\gvac{1}\gob{1}{A}
\gend
=\gbeg{1}{3}
\got{1}{A}\gnl
\gcl{1}\gnl
\gob{1}{A}
\gend
$ 
we find that $\un{\Delta}_A$ respects the unit of $A$, that is, the sixth equality in
\equref{unitcounitrel} is satisfied. The second equality in \equref{unitcondtensprod} is
equivalent to 
\[
\gbeg{3}{7}
\gvac{1}\got{1}{Y}\gnl
\gvac{1}\gcl{2}\gnl
\gdb\gnl
\gcl{1}\gcn{1}{1}{1}{3}\gnl
\gcl{1}\gu{1}\gcl{1}\gnl
\glm\gcl{1}\gnl
\gvac{1}\gob{1}{A}\gob{1}{Y}
\gend
=
\gbeg{2}{3}
\gvac{1}\got{1}{Y}\gnl
\gu{1}\gcl{1}\gnl
\gob{1}{A}\gob{1}{Y}
\gend
\]
and this is equivalent to the last equality in \equref{unitcounitrel}; for the direct
implication, take $Y= B_{\rm tr}$, and compose to the left ${\rm id}_A\ot \un{\va}_B$.
This finishes our proof.
\end{proof}

Our next aim is to show that the left $B$-action on $A$ defines a left $B$-module algebra
and a left $B$-module coalgebra structure on $A$. First we remark that the right $A$-module
structure on $X\ot Y$ ($X,Y\in {}^B\Cc_A$) satisfies the associativity condition if and only if
\[
\gbeg{5}{9}
\got{1}{X}\gvac{1}\got{1}{Y}\got{1}{A}\got{1}{A}\gnl
\gcl{1}\gvac{1}\gcl{2}\gmu\gnl
\gcl{1}\gdb\gcmu\gnl
\gcl{1}\gcl{1}\gbr\gcl{1}\gnl
\gcl{1}\glm\gcl{1}\gvac{-1}\gev\gnl
\gcl{1}\gcn{1}{1}{3}{1}\gvac{1}\gcl{1}\gnl
\gcl{2}\gvac{-1}\gev\gvac{1}\gcl{2}\gnl
\gob{1}{X}\gvac{2}\gob{1}{Y}
\gend
=
\gbeg{6}{13}
\got{1}{X}\gvac{1}\got{1}{Y}\got{2}{A}\got{0}{A}\gnl
\gcl{1}\gvac{1}\gcl{2}\gcmu\gcn{1}{4}{0}{0}\gnl
\gcl{1}\gdb\gcl{1}\gcl{1}\gnl
\gcl{1}\gcl{1}\gbr\gcl{1}\gnl
\gcl{1}\glm\gcl{2}\gvac{-1}\gev\gnl
\gcl{1}\gcn{1}{1}{3}{1}\gvac{2}\gcmu\gnl
\gcl{2}\gvac{-1}\gev\gdb\gvac{-1}\gcl{1}\gcl{1}\gcl{1}\gnl
\gvac{2}\gcl{1}\gbr\gcl{1}\gnl
\gcl{1}\gvac{1}\glm\gcl{2}\gvac{-1}\gev\gnl
\gcl{1}\gvac{2}\gcn{1}{1}{1}{-3}\gnl
\gcl{2}\gvac{-1}\gev\gvac{2}\gcl{2}\gnl
\gob{1}{X}\gvac{3}\gob{1}{Y}
\gend\hspace{1mm}.
\]
Since $X$ is a right $A$-module, this is equivalent to
\begin{equation}\eqlabel{rAmodcondtens}
\gbeg{4}{7}
\gvac{1}\got{1}{Y}\got{1}{A}\got{1}{A}\gnl
\gvac{1}\gcl{2}\gmu\gnl
\gdb\gcmu\gnl
\gcl{1}\gbr\gcl{1}\gnl
\glm\gcl{2}\gvac{-1}\gev\gnl
\gvac{1}\gcl{1}\gnl
\gvac{1}\gob{1}{A}\gob{1}{Y}
\gend
=
\gbeg{5}{13}
\gvac{1}\got{1}{Y}\got{2}{A}\got{0}{A}\gnl
\gvac{1}\gcl{2}\gcmu\gcn{1}{4}{0}{0}\gnl
\gdb\gcl{1}\gcl{1}\gnl
\gcl{1}\gbr\gcl{1}\gnl
\glm\gcl{2}\gvac{-1}\gev\gnl
\gvac{1}\gcn{1}{1}{1}{-1}\gvac{1}\gcmu\gnl
\gcl{1}\gdb\gvac{-1}\gcl{2}\gcl{1}\gcl{1}\gnl
\gcl{1}\gcl{1}\gcl{1}\gcl{1}\gcl{1}\gnl
\gcl{1}\gcl{1}\gbr\gcl{1}\gnl
\gcl{1}\glm\gcl{3}\gvac{-1}\gev\gnl
\gcl{1}\gcn{1}{1}{3}{1}\gnl
\gmu\gnl
\gob{2}{A}\gvac{1}\gob{1}{Y}
\gend
\hspace{2mm},\mbox{$\forall~~Y\in {}^B\Cc_A$}\hspace{1mm}.
\end{equation}
For the direct implication, take $X=A$, and compose to the right with $\un{\eta}_A\ot
{\rm id}_{Y\ot A\ot A}$.\\
In a similar way, it can be shown that the associativity constraints $a_{X,Y,Z}$ are
right $A$-linear if and only if
\begin{equation}\eqlabel{rAmodcondass}
\gbeg{6}{11}
\gvac{1}\got{1}{Y}\gvac{1}\got{1}{Z}\got{2}{A}\gnl
\gvac{1}\gcl{2}\gvac{1}\gcl{2}\gcmu\gnl
\gdb\gdb\gcl{1}\gcl{1}\gnl
\gcl{1}\gcl{1}\gcl{1}\gbr\gcl{1}\gnl
\gcl{1}\gcl{1}\glm\gcl{6}\gvac{-1}\gev\gnl
\gcl{1}\gcl{1}\gvac{1}\gcn{1}{1}{1}{0}\gnl
\gcl{1}\gcl{1}\gcmu\gnl
\gcl{1}\gbr\gcl{1}\gnl
\glm\gcl{2}\gvac{-1}\gev\gnl
\gvac{1}\gcl{1}\gnl
\gvac{1}\gob{1}{A}\gob{1}{Y}\gvac{1}\gob{1}{Z}
\gend
=
\gbeg{7}{14}
\gvac{1}\got{1}{Y}\gvac{1}\got{1}{Z}\got{2}{A}\gnl
\gvac{1}\gcl{2}\gvac{1}\gcl{2}\gcmu\gnl
\gdb\gdb\gcl{1}\gcl{1}\gnl
\gcl{1}\gbr\gbr\gcn{1}{1}{1}{2}\gnl
\gmu\gbr\gcl{2}\gcmu\gnl
\gcn{1}{1}{2}{1}\gvac{1}\gcn{1}{1}{1}{-1}\gcn{1}{1}{1}{-1}\gvac{-1}\gdb\gvac{-1}\gcl{1}\gcl{1}\gcl{1}\gnl
\glm\gcl{1}\gcl{1}\gcl{1}\gcl{1}\gcl{1}\gnl
\gvac{1}\gcl{1}\gcl{1}\gcl{1}\gbr\gcl{1}\gnl
\gvac{1}\gcl{1}\gcl{1}\glm\gcl{2}\gvac{-1}\gev\gnl
\gvac{1}\gcl{1}\gcl{1}\gvac{1}\gcl{1}\gnl
\gvac{1}\gcl{1}\gcl{1}\gcn{1}{1}{3}{1}\gvac{1}\gcl{2}\gnl
\gvac{1}\gcl{1}\gcl{2}\gvac{-1}\gev\gnl
\gvac{1}\gcl{1}\gvac{3}\gcl{1}\gnl
\gvac{1}\gob{1}{A}\gob{1}{Y}\gvac{2}\gob{1}{Z}
\gend,
\end{equation}
for all $Y, Z\in {}^B\Cc_A$. 
The verification is left to the reader.

\begin{lemma}\lelabel{2.4}
Assume that the input monoidal structure $({}^B\Cc_A, \ot, \un{1}, a, l, r)$
associated to an input monoidal Hopf module datum $(B, A)$ is monoidal.
Then the left $B$-action on $A$ makes $A$ into a left $B$-module algebra
and a left $B$-module coalgebra.
\end{lemma}

\begin{proof}
We first show that $A$ is a left $B$-module. Taking $Y=Z=B_{\rm tr}$ in \equref{rAmodcondass},
we obtain that
\[
\gbeg{5}{8}
\got{2}{B}\got{2}{B}\got{1}{A}\gnl
\gcmu\gcmu\gcl{1}\gnl
\gcl{1}\gcl{1}\gcl{1}\gbr\gnl
\gcl{1}\gcl{1}\glm\gcl{1}\gnl
\gcl{1}\gcl{1}\gvac{1}\gcn{1}{1}{1}{-1}\gcl{1}\gnl
\gcl{1}\gbr\gvac{1}\gcl{1}\gnl
\glm\gcl{1}\gvac{1}\gcl{1}\gnl
\gvac{1}\gob{1}{A}\gob{1}{B}\gvac{1}\gob{1}{B}
\gend
=
\gbeg{5}{7}
\got{2}{B}\got{2}{B}\got{1}{A}\gnl
\gcmu\gcmu\gcl{1}\gnl
\gcl{1}\gbr\gbr\gnl
\gmu\gbr\gcl{1}\gnl
\gcn{1}{1}{2}{1}\gvac{1}\gcn{1}{1}{1}{-1}\gcl{1}\gcl{1}\gnl
\glm\gvac{1}\gcl{1}\gcl{1}\gnl
\gvac{1}\gob{1}{A}\gvac{1}\gob{1}{B}\gob{1}{B}
\gend
\hspace{1mm}.
\] 
Composing this identity to the left with ${\rm id}_A\ot \un{\va}_B\ot \un{\va}_B$, we find that
$$
\gbeg{3}{5}
\got{1}{B}\got{1}{B}\got{1}{A}\gnl
\gcl{1}\glm\gnl
\gcl{1}\gcn{1}{1}{3}{1}\gnl
\glm\gnl
\gvac{1}\gob{1}{A}
\gend
=
\gbeg{3}{5}
\got{1}{B}\got{1}{B}\got{1}{A}\gnl
\gmu\gcl{1}\gnl
\gcn{1}{1}{2}{3}\gvac{1}\gcl{1}\gnl
\gvac{1}\glm\gnl
\gvac{2}\gob{1}{A}
\gend\hspace{1mm}.
$$
Together with the fifth equality in \equref{unitcounitrel}, this shows that $A$ is a left $B$-module.\\
Now take $Y=B_{\rm tr}$ in \equref{rAmodcondtens}, and compose at the left with
$\un{\va}_B\ot {\rm id}_A$, to obtain the second equality in \equref{lmodalg}.
Together with the last equality in \equref{unitcounitrel}, this tells us that $A$ is a left $B$-module 
algebra.\\
Now take $Y=A$ and $Z=B_{\rm tr}$ in \equref{rAmodcondass}, and compose to the left
with $ {\rm id}_A\ot  {\rm id}_A\ot \un{\va}_B$. Then we obtain
\[
\gbeg{5}{7}
\gvac{1}\got{1}{A}\got{1}{B}\got{1}{A}\gnl
\glcm\glm\gnl
\gcl{1}\gcl{1}\gcn{1}{1}{3}{2}\gnl
\gcl{1}\gcl{1}\gcmu\gnl
\gcl{1}\gbr\gcl{1}\gnl
\glm\gmu\gnl
\gvac{1}\gob{1}{A}\gob{2}{A}
\gend
=
\gbeg{6}{7}
\gvac{1}\got{1}{A}\got{2}{B}\got{2}{A}\gnl
\glcm\gcmu\gcmu\gnl
\gcl{1}\gbr\gbr\gcl{1}\gnl
\gmu\gbr\glm\gnl
\gcn{1}{1}{2}{1}\gvac{1}\gcn{1}{1}{1}{-1}\gcl{1}\gcn{1}{1}{3}{1}\gnl
\glm\gvac{1}\gmu\gnl
\gvac{1}\gob{1}{A}\gvac{1}\gob{2}{A}
\gend
\hspace{1mm}.
\]
Now compose to the right with $\un{\eta}_A\ot {\rm id}_B\ot {\rm id}_A$. Taking into account that
$
\gbeg{2}{4}
\got{2}{\un{1}}\gnl
\gvac{1}\gu{1}\gnl
\glcm\gnl
\gob{1}{B}\gob{1}{A}
\gend
=
\gbeg{2}{3}
\got{2}{\un{1}}\gnl
\gu{1}\gu{1}\gnl
\gob{1}{B}\gob{1}{A}
\gend
$, 
we find the second equality in  \equref{lmodcoalg}. Combined with the fourth equality
in \equref{unitcounitrel}, this tells us that $A$ is a left $B$-module coalgebra.
\end{proof}

Our next step is to show that $A$ is a bialgebra in ${}^B_B{\cal YD}$. To this end, we
first point out that the tensor product $X\ot Y$ of two objects $X, Y\in {}^B\Cc_A$
satisfies the compatibility relation for relative Hopf modules if and only if
\begin{equation}\eqlabel{tensprodcomprel}
\gbeg{6}{13}
\gvac{1}\got{1}{X}\gvac{1}\got{1}{Y}\got{2}{A}\gnl
\gvac{1}\gcl{1}\gvac{1}\gcl{2}\gcmu\gnl
\gvac{1}\gcl{1}\gdb\gcl{1}\gcl{1}\gnl
\gvac{1}\gcl{1}\gcl{1}\gbr\gcl{1}\gnl
\gvac{1}\gcl{1}\glm\gcl{2}\gvac{-1}\gev\gnl
\gvac{1}\gcl{1}\gcn{1}{1}{3}{1}\gvac{1}\gnl
\gvac{1}\gcl{2}\gvac{-1}\gev\gdb\gvac{-1}\gcl{2}\gnl
\gvac{1}\gcl{2}\gvac{1}\gcl{1}\gnl
\gdb\gvac{-1}\gcl{2}\gvac{1}\gcn{1}{2}{1}{-1}\gcl{4}\gnl
\gcl{1}\gnl
\gcl{1}\gbr\gnl
\gmu\gcl{1}\gnl
\gob{2}{B}\gob{1}{X}\gvac{1}\gob{1}{Y}
\gend
=
\gbeg{8}{13}
\gvac{1}\got{1}{X}\gvac{1}\got{1}{Y}\gvac{2}\got{1}{A}\gnl
\gvac{1}\gcl{2}\gvac{1}\gcl{2}\gvac{1}\glcm\gnl
\gdb\gdb\gvac{1}\gcl{1}\gcl{1}\gnl
\gcl{1}\gbr\gcn{1}{1}{1}{3}\gvac{1}\gcl{1}\gcl{1}\gnl
\gmu\gcl{1}\gvac{1}\gbr\gcn{1}{1}{1}{2}\gnl
\gcn{1}{2}{2}{3}\gvac{1}\gcl{1}\gvac{1}\gcn{1}{1}{1}{-1}\gcl{2}\gcmu\gnl
\gvac{2}\gbr\gdb\gcl{1}\gcl{1}\gnl
\gvac{1}\gmu\gcl{1}\gcl{1}\gbr\gcl{1}\gnl
\gvac{1}\gcn{1}{4}{2}{2}\gvac{1}\gcl{1}\glm\gcl{4}\gvac{-1}\gev\gnl
\gvac{3}\gcl{1}\gvac{1}\gcn{1}{1}{1}{-1}\gnl
\gvac{3}\gcl{2}\gvac{-1}\gev\gnl
\gvac{1}\gob{2}{B}\gob{1}{X}\gvac{2}\gob{1}{Y}
\gend
\hspace{1mm}.
\end{equation} 

\begin{proposition}\prlabel{directimplic}
If the input monoidal structure $({}^B\Cc_A, \ot, \un{1}, a, l, r)$
associated to an input monoidal Hopf module datum $(B, A)$ is monoidal,
then $A$ is a bialgebra in ${}^B_B{\cal YD}$.
\end{proposition}

\begin{proof}
It follows from our previous results that it suffices to show that
$A$ is a left $B$-comodule coalgebra and an object of ${}^B_B{\cal YD}$,
and that $\un{\Delta}_A$ satisfies the last equality in \equref{bialgYD}.\\
We take $X=Y=A$ in \equref{tensprodcomprel}, and then compose to the right
with $\un{\eta}_A\ot \un{\eta}_A\ot {\rm id}_A$. Using the first relation in 
\equref{lcomalg} and the fact that $A$ is a left $B$-module, we obtain
the second equality in \equref{lcomcoalg}. The first equality in \equref{lcomcoalg}
was proved in \leref{unitcounitcomp}, hence it follows that 
$A$ is a left $B$-comodule coalgebra.\\
Now let $X=A$ and $Y=B_{\rm tr}$ in \equref{tensprodcomprel}, and compose to the
right with $\un{\eta}_A\ot {\rm id}_B \ot {\rm id}_A$, and to the left with ${\rm id}_B \ot {\rm id}_A
\ot \un{\va}_B$. Then we obtain that
\begin{equation}\eqlabel{YDcondForB}
\gbeg{3}{8}
\got{2}{B}\got{1}{A}\gnl
\gcmu\gcl{1}\gnl
\gcl{1}\gbr\gnl
\glm\gcl{1}\gnl
\glcm\gcl{1}\gnl
\gcl{1}\gbr\gnl
\gmu\gcl{1}\gnl
\gob{2}{B}\gob{1}{A}
\gend
=
\gbeg{4}{5}
\got{2}{B}\gvac{1}\got{1}{A}\gnl
\gcmu\glcm\gnl
\gcl{1}\gbr\gcl{1}\gnl
\gmu\glm\gnl
\gob{2}{B}\gvac{1}\gob{1}{A}
\gend
\hspace{1mm},
\end{equation} 
which is precisely the required compatibility condition between the left $B$-action and left
$B$-coaction on $A$ that is needed to make $A$ a left Yetter-Drinfeld module over $B$ in $\Cc$.\\
Finally, take $Y=A$ in \equref{rAmodcondtens}, and compose to the right with
$\un{\eta}_A\ot {\rm id}_A\ot {\rm id}_A$. 
Due to 
the first equality in \equref{lcomalg} and the fact that $A$ is a left $B$-module, it turns out that the last 
equality in \equref{bialgYD} holds, and this finishes our proof.
\end{proof}

{\sl Proof of \thref{2.1}}.
One implication has been proved in \prref{directimplic}. Conversely, assume that
$A$ is a bialgebra in ${}^B_B{\cal YD}$. Since ${}^B\Cc$ is a monoidal category ($B$
is a bialgebra in $\Cc$), we have that $l_X$, $r_X$ and $a_{X,Y,Z}$ are left
$B$-colinear, for all $X,Y,X\in {}^B\Cc_A$. From \leref{unitcounitcomp}, we know
that $\un{1}\in {}^B\Cc_A$, that $l_X$ and $r_X$ are also right $A$-linear, and 
the right $A$-action on $X\ot Y$ defined in \equref{rAactDHM} is unital. Thus it remains to be shown that
$a_{X,Y,Z}$ is right $A$-linear, that $X\ot Y$ is associative as a right $A$-module,
and that it satisfies the compatibilty condition for a relative Hopf module. Otherwise
stated, we have to show that $\ref{eq:rAmodcondtens}-\ref{eq:tensprodcomprel}$
are satisfied.\\
To prove \equref{rAmodcondtens} we compute, for all $Y\in {}^B\Cc_A$, 
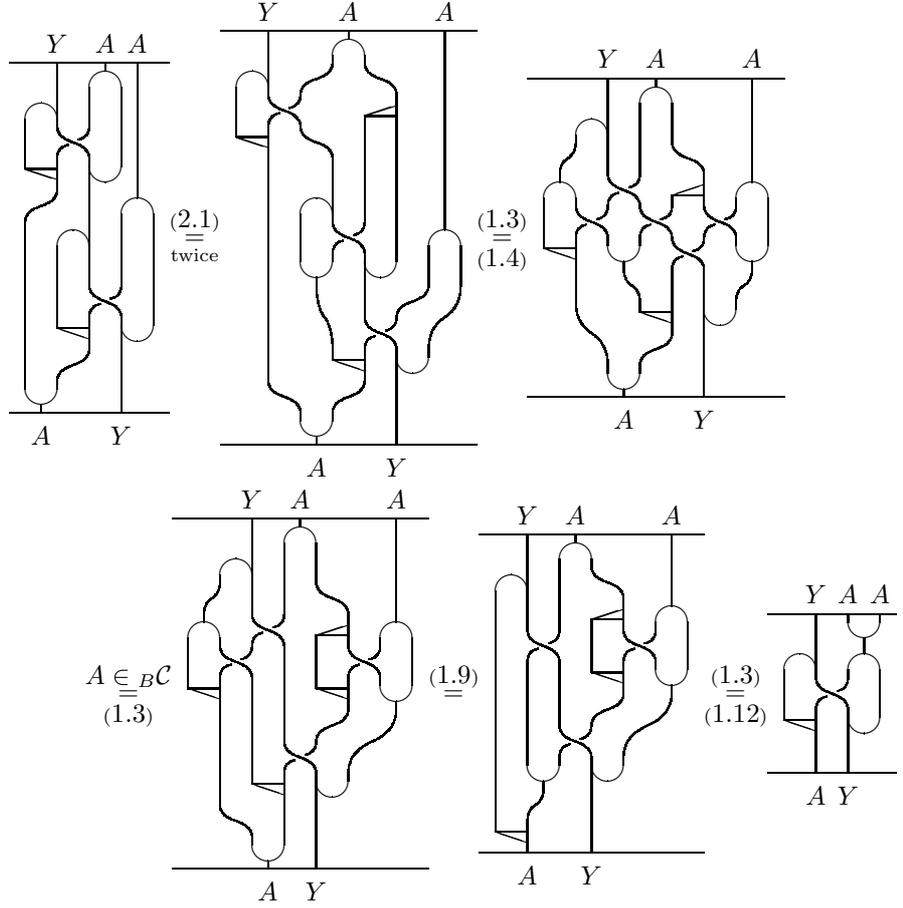
\begin{eqnarray*}
&&\hspace*{-1cm}
\gbeg{5}{13}
\gvac{1}\got{1}{Y}\got{2}{A}\got{0}{A}\gnl
\gvac{1}\gcl{2}\gcmu\gcn{1}{4}{0}{0}\gnl
\gdb\gcl{1}\gcl{1}\gnl
\gcl{1}\gbr\gcl{1}\gnl
\glm\gcl{2}\gvac{-1}\gev\gnl
\gvac{1}\gcn{1}{1}{1}{-1}\gvac{1}\gcmu\gnl
\gcl{1}\gdb\gvac{-1}\gcl{2}\gcl{1}\gcl{1}\gnl
\gcl{1}\gcl{1}\gcl{1}\gcl{1}\gcl{1}\gnl
\gcl{1}\gcl{1}\gbr\gcl{1}\gnl
\gcl{1}\glm\gcl{3}\gvac{-1}\gev\gnl
\gcl{1}\gcn{1}{1}{3}{1}\gnl
\gmu\gnl
\gob{2}{A}\gvac{1}\gob{1}{Y}
\gend
\equalupdown{\equref{rDHcond}}{\rm twice}
\gbeg{8}{15}
\gvac{1}\got{1}{Y}\gvac{1}\got{2}{A}\gvac{1}\got{2}{A}\gnl
\gvac{1}\gcl{2}\gvac{1}\gcmu\gvac{1}\gcn{1}{6}{2}{2}\gnl
\gdb\gvac{1}\gcn{1}{1}{1}{-1}\gcn{1}{1}{1}{3}\gnl
\gcl{1}\gbr\gvac{1}\glcm\gnl
\glm\gcn{1}{1}{1}{3}\gvac{1}\gcl{1}\gcl{1}\gnl
\gvac{1}\gcl{1}\gvac{1}\gcl{2}\gcl{1}\gcl{1}\gnl
\gvac{1}\gcl{1}\gdb\gcl{1}\gcl{1}\gnl
\gvac{1}\gcl{1}\gcl{1}\gbr\gcl{1}\gcmu\gnl
\gvac{1}\gcl{1}\gmu\gcl{2}\gvac{-1}\gev\gcl{1}\gcl{1}\gnl
\gvac{1}\gcl{3}\gcn{1}{2}{2}{3}\gvac{3}\gcn{1}{1}{1}{-1}\gcn{1}{2}{1}{-1}\gnl
\gvac{4}\gbr\gnl
\gvac{3}\glm\gcl{3}\gvac{-1}\gev\gnl
\gvac{1}\gcn{1}{1}{1}{3}\gvac{1}\gcn{1}{1}{3}{1}\gnl
\gvac{2}\gmu\gnl
\gvac{2}\gob{2}{A}\gvac{1}\gob{1}{Y}
\gend
\equalupdown{\equref{nat1cup}}{\equref{nat2cup}}
\gbeg{8}{12}
\gvac{2}\got{1}{Y}\got{2}{A}\gvac{1}\got{2}{A}\gnl
\gvac{2}\gcl{2}\gcmu\gvac{1}\gcn{1}{3}{2}{2}\gnl
\gvac{1}\gdb\gcl{1}\gcl{1}\gnl
\gvac{1}\gcn{1}{1}{1}{0}\gcl{1}\gcl{1}\gcn{1}{1}{1}{3}\gnl
\gcmu\gbr\glcm\gcmu\gnl
\gcl{1}\gbr\gbr\gbr\gcl{1}\gnl
\glm\gmu\gbr\gmu\gnl
\gvac{1}\gcl{1}\gcn{1}{1}{2}{3}\gvac{1}\gcl{1}\gcl{2}\gcn{1}{1}{2}{1}\gnl
\gvac{1}\gcn{1}{2}{1}{3}\gvac{1}\glm\gev\gnl
\gvac{4}\gcn{1}{1}{1}{-1}\gcl{2}\gnl
\gvac{2}\gmu\gnl
\gvac{2}\gob{2}{A}\gvac{1}\gob{1}{Y}
\gend\\
&&\equalupdown{\mbox{$A\in {}_B\Cc$}}{\equref{nat1cup}}
\gbeg{8}{13}
\gvac{2}\got{1}{Y}\got{2}{A}\gvac{1}\got{2}{A}\gnl
\gvac{2}\gcl{2}\gcmu\gvac{1}\gcn{1}{3}{2}{2}\gnl
\gvac{1}\gdb\gcl{1}\gcl{1}\gnl
\gvac{1}\gcn{1}{1}{1}{0}\gcl{1}\gcl{1}\gcn{1}{1}{1}{3}\gnl
\gcmu\gbr\glcm\gcmu\gnl
\gcl{1}\gbr\gcl{1}\gcl{1}\gbr\gcl{1}\gnl
\glm\gcl{1}\gcl{1}\glm\gmu\gnl
\gvac{1}\gcl{1}\gcl{1}\gcl{1}\gcn{1}{1}{3}{1}\gvac{1}\gcn{1}{2}{2}{-1}\gnl
\gvac{1}\gcl{1}\gcl{1}\gbr\gnl
\gvac{1}\gcl{1}\glm\gcl{3}\gvac{-1}\gev\gnl
\gvac{1}\gcn{1}{1}{1}{3}\gvac{1}\gcl{1}\gnl
\gvac{2}\gmu\gnl
\gvac{2}\gob{2}{A}\gob{1}{Y}
\gend
\equal{\equref{lmodalg}}
\gbeg{7}{12}
\gvac{1}\got{1}{Y}\got{2}{A}\gvac{1}\got{2}{A}\gnl
\gvac{1}\gcl{2}\gcmu\gvac{1}\gcn{1}{2}{2}{2}\gnl
\gdb\gcl{1}\gcn{1}{1}{1}{3}\gnl
\gcl{1}\gcl{1}\gcl{1}\glcm\gcmu\gnl
\gcl{1}\gbr\gcl{1}\gbr\gcl{1}\gnl
\gcl{1}\gcl{1}\gcl{1}\glm\gmu\gnl
\gcl{1}\gcl{1}\gcl{1}\gcn{1}{1}{3}{1}\gvac{1}\gcn{1}{2}{2}{-1}\gnl
\gcl{1}\gcl{1}\gbr\gnl
\gcl{1}\gmu\gcl{3}\gvac{-1}\gev\gnl
\gcl{1}\gcn{1}{1}{2}{1}\gnl
\glm\gnl
\gvac{1}\gob{1}{A}\gvac{1}\gob{1}{Y}
\gend
\equalupdown{\equref{nat1cup}}{\equref{bialgYD}}
\gbeg{4}{7}
\gvac{1}\got{1}{Y}\got{1}{A}\got{1}{A}\gnl
\gvac{1}\gcl{2}\gmu\gnl
\gdb\gcmu\gnl
\gcl{1}\gbr\gcl{1}\gnl
\glm\gcl{2}\gvac{-1}\gev\gnl
\gvac{1}\gcl{1}\gnl
\gvac{1}\gob{1}{A}\gob{1}{Y}
\gend
\hspace{1mm},
\end{eqnarray*}
as needed. The proof of \equref{rAmodcondass} is similar, and is left to the reader.
Observe that is essentially based on the fact that $A$ is a left $B$-module coalgebra.\\
Finally, for all $X, Y\in {}^B\Cc_A$, we have that
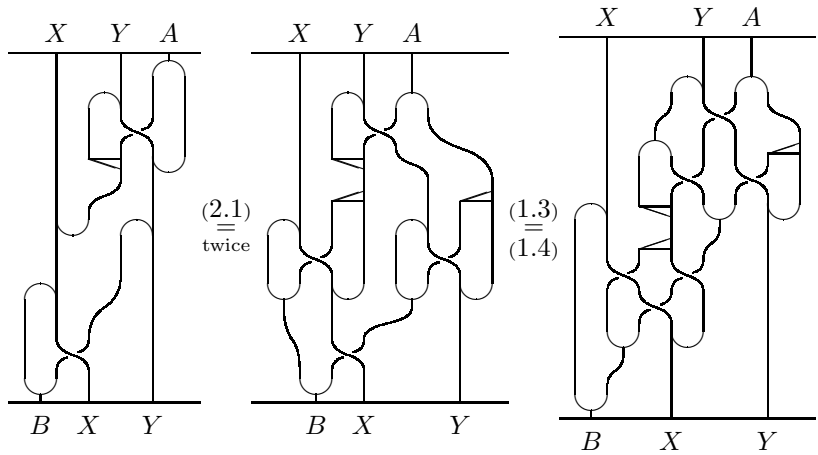
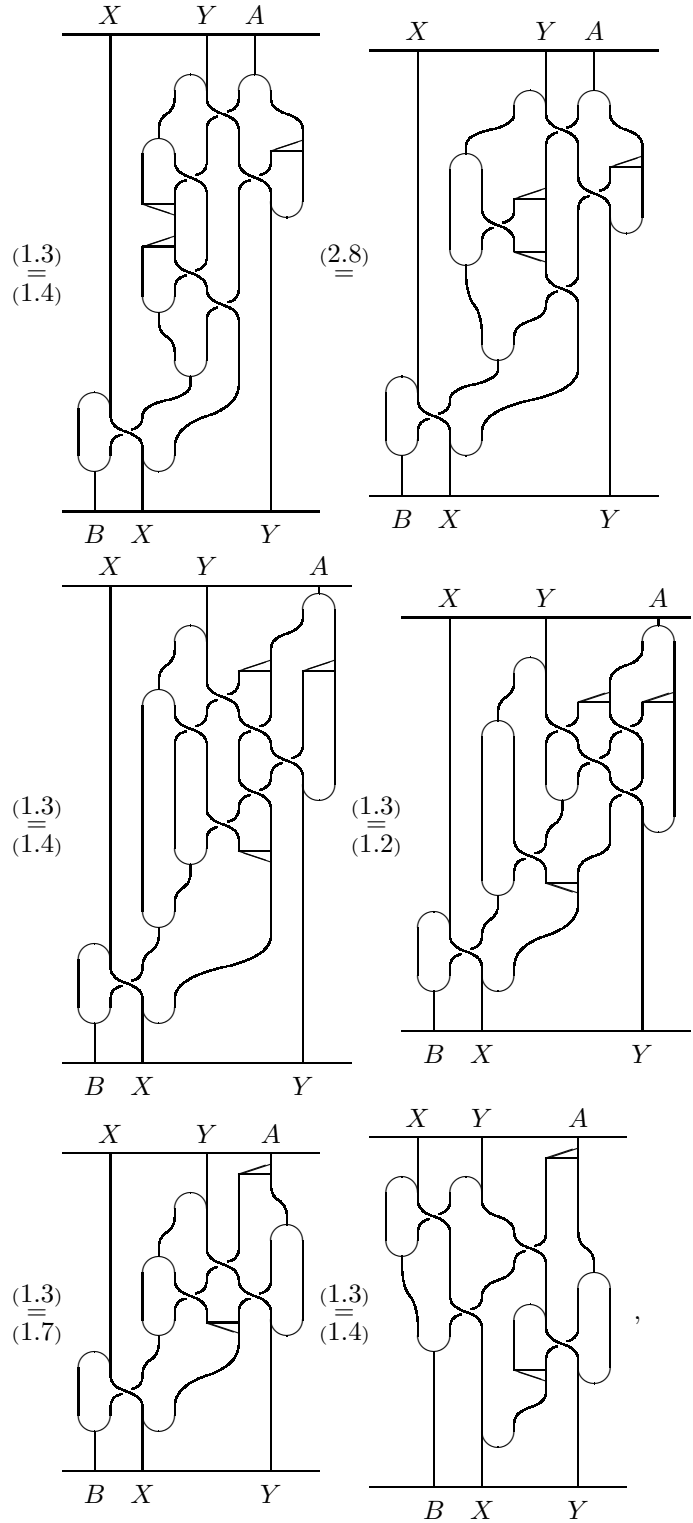
\begin{eqnarray*}
&&\hspace*{-1cm}
\gbeg{6}{13}
\gvac{1}\got{1}{X}\gvac{1}\got{1}{Y}\got{2}{A}\gnl
\gvac{1}\gcl{1}\gvac{1}\gcl{2}\gcmu\gnl
\gvac{1}\gcl{1}\gdb\gcl{1}\gcl{1}\gnl
\gvac{1}\gcl{1}\gcl{1}\gbr\gcl{1}\gnl
\gvac{1}\gcl{1}\glm\gcl{2}\gvac{-1}\gev\gnl
\gvac{1}\gcl{1}\gcn{1}{1}{3}{1}\gvac{1}\gnl
\gvac{1}\gcl{2}\gvac{-1}\gev\gdb\gvac{-1}\gcl{2}\gnl
\gvac{1}\gcl{2}\gvac{1}\gcl{1}\gnl
\gdb\gvac{-1}\gcl{2}\gvac{1}\gcn{1}{2}{1}{-1}\gcl{4}\gnl
\gcl{1}\gnl
\gcl{1}\gbr\gnl
\gmu\gcl{1}\gnl
\gob{2}{B}\gob{1}{X}\gvac{1}\gob{1}{Y}
\gend
\equalupdown{\equref{rDHcond}}{\rm twice}
\gbeg{8}{13}
\gvac{1}\got{1}{X}\gvac{1}\got{1}{Y}\got{2}{A}\gnl
\gvac{1}\gcl{1}\gvac{1}\gcl{2}\gcn{2}{1}{2}{2}\gnl
\gvac{1}\gcl{1}\gdb\gcmu\gnl
\gvac{1}\gcl{1}\gcl{1}\gbr\gcn{1}{2}{1}{5}\gnl
\gvac{1}\gcl{1}\glm\gcn{1}{1}{1}{3}\gnl
\gvac{1}\gcl{2}\glcm\gvac{1}\gcl{2}\glcm\gnl
\gdb\gcl{1}\gcl{1}\gdb\gcl{1}\gcl{1}\gnl
\gcl{1}\gbr\gcl{1}\gcl{1}\gbr\gcl{1}\gnl
\gmu\gcl{2}\gvac{-1}\gev\gmu\gcl{4}\gvac{-1}\gev\gnl
\gcn{1}{2}{2}{3}\gvac{2}\gcn{1}{1}{4}{1}\gnl
\gvac{2}\gbr\gnl
\gvac{1}\gmu\gcl{1}\gnl
\gvac{1}\gob{2}{B}\gob{1}{X}\gvac{2}\gob{1}{Y}
\gend
\equalupdown{\equref{nat1cup}}{\equref{nat2cup}}
\gbeg{8}{14}
\gvac{1}\got{1}{X}\gvac{2}\got{1}{Y}\got{2}{A}\gnl
\gvac{1}\gcl{1}\gvac{2}\gcl{2}\gcn{2}{1}{2}{2}\gnl
\gvac{1}\gcl{1}\gvac{1}\gdb\gcmu\gnl
\gvac{1}\gcl{1}\gvac{1}\gcn{1}{1}{1}{0}\gbr\gcn{1}{1}{1}{3}\gnl
\gvac{1}\gcl{1}\gcmu\gcl{1}\gcl{1}\glcm\gnl
\gvac{1}\gcl{1}\gcl{1}\gbr\gbr\gcl{1}\gnl
\gdb\gvac{-1}\gcl{1}\glm\gmu\gcl{7}\gvac{-1}\gev\gnl
\gcl{1}\gcl{1}\glcm\gcn{1}{1}{2}{1}\gnl
\gcl{1}\gbr\gbr\gnl
\gcl{1}\gcl{1}\gbr\gcl{1}\gnl
\gcl{1}\gmu\gcl{3}\gvac{-1}\gev\gnl
\gcl{1}\gcn{1}{1}{2}{1}\gnl
\gmu\gnl
\gob{2}{B}\gvac{1}\gob{1}{X}\gvac{2}\gob{1}{Y}
\gend\\
&&
\equalupdown{\equref{nat1cup}}{\equref{nat2cup}}
\gbeg{8}{17}
\gvac{1}\got{1}{X}\gvac{2}\got{1}{Y}\got{2}{A}\gnl
\gvac{1}\gcl{1}\gvac{2}\gcl{2}\gcn{2}{1}{2}{2}\gnl
\gvac{1}\gcl{1}\gvac{1}\gdb\gcmu\gnl
\gvac{1}\gcl{1}\gvac{1}\gcn{1}{1}{1}{0}\gbr\gcn{1}{1}{1}{3}\gnl
\gvac{1}\gcl{1}\gcmu\gcl{1}\gcl{1}\glcm\gnl
\gvac{1}\gcl{1}\gcl{1}\gbr\gbr\gcl{1}\gnl
\gvac{1}\gcl{1}\glm\gcl{1}\gcl{1}\gcl{10}\gvac{-1}\gev\gnl
\gvac{1}\gcl{1}\glcm\gcl{1}\gcl{1}\gnl
\gvac{1}\gcl{1}\gcl{1}\gbr\gcl{1}\gnl
\gvac{1}\gcl{1}\gmu\gbr\gnl
\gvac{1}\gcl{1}\gcn{1}{1}{2}{3}\gvac{1}\gcl{1}\gcl{1}\gnl
\gvac{1}\gcl{1}\gvac{1}\gmu\gcl{1}\gnl
\gdb\gvac{-1}\gcl{1}\gvac{2}\gcn{1}{1}{0}{-3}\gcn{1}{2}{1}{-3}\gnl
\gcl{1}\gbr\gnl
\gmu\gcl{2}\gvac{-1}\gev\gnl
\gcn{1}{1}{2}{2}\gnl
\gob{2}{B}\gob{1}{X}\gvac{3}\gob{1}{Y}
\gend
\equal{\equref{YDcondForB}}
\gbeg{9}{16}
\gvac{1}\got{1}{X}\gvac{3}\got{1}{Y}\got{2}{A}\gnl
\gvac{1}\gcl{1}\gvac{3}\gcl{2}\gcn{2}{1}{2}{2}\gnl
\gvac{1}\gcl{1}\gvac{2}\gdb\gcmu\gnl
\gvac{1}\gcl{1}\gvac{2}\gcn{1}{1}{1}{-2}\gbr\gcn{1}{1}{1}{3}\gnl
\gvac{1}\gcl{1}\gcmu\gvac{1}\gcl{1}\gcl{1}\glcm\gnl
\gvac{1}\gcl{1}\gcl{1}\gcl{1}\glcm\gbr\gcl{1}\gnl
\gvac{1}\gcl{1}\gcl{1}\gbr\gcl{1}\gcl{1}\gcl{9}\gvac{-1}\gev\gnl
\gvac{1}\gcl{1}\gmu\glm\gcl{1}\gnl
\gvac{1}\gcl{1}\gcn{1}{2}{2}{3}\gvac{2}\gbr\gnl
\gvac{1}\gcl{1}\gvac{3}\gcn{1}{1}{1}{-1}\gcl{2}\gnl
\gvac{1}\gcl{1}\gvac{1}\gmu\gnl
\gdb\gvac{-1}\gcl{1}\gvac{1}\gcn{1}{1}{2}{-1}\gvac{2}\gcn{1}{2}{1}{-5}\gnl
\gcl{1}\gbr\gnl
\gmu\gcl{2}\gvac{-1}\gev\gnl
\gcn{1}{1}{2}{2}\gnl
\gob{2}{B}\gob{1}{X}\gvac{4}\gob{1}{Y}
\gend\\
&&\equalupdown{\equref{nat1cup}}{\equref{nat2cup}}
\gbeg{9}{17}
\gvac{1}\got{1}{X}\gvac{2}\got{1}{Y}\gvac{2}\got{2}{A}\gnl
\gvac{1}\gcl{1}\gvac{2}\gcl{2}\gvac{2}\gcmu\gnl
\gvac{1}\gcl{1}\gvac{1}\gdb\gvac{2}\gcn{1}{1}{1}{-1}\gcl{1}\gnl
\gvac{1}\gcl{1}\gvac{1}\gcn{1}{1}{1}{0}\gcl{1}\glcm\glcm\gnl
\gvac{1}\gcl{1}\gcmu\gbr\gcl{1}\gcl{1}\gcl{1}\gnl
\gvac{1}\gcl{1}\gcl{1}\gbr\gbr\gcl{1}\gcl{1}\gnl
\gvac{1}\gcl{1}\gcl{5}\gcl{2}\gcl{1}\gcl{1}\gbr\gcl{1}\gnl
\gvac{1}\gcl{1}\gvac{2}\gcl{1}\gbr\gcl{9}\gvac{-1}\gev\gnl
\gvac{1}\gcl{1}\gvac{1}\gcl{1}\gbr\gcl{1}\gnl
\gvac{1}\gcl{1}\gvac{1}\gmu\glm\gnl
\gvac{1}\gcl{1}\gvac{1}\gcn{1}{1}{2}{1}\gvac{2}\gcl{2}\gnl
\gvac{1}\gcl{1}\gmu\gnl
\gdb\gvac{-1}\gcl{1}\gcn{1}{1}{2}{1}\gvac{2}\gcn{1}{2}{3}{-3}\gnl
\gcl{1}\gbr\gnl
\gmu\gcl{2}\gvac{-1}\gev\gnl
\gcn{1}{1}{2}{2}\gnl
\gob{2}{B}\gob{1}{X}\gvac{4}\gob{1}{Y}
\gend
\equalupdown{\equref{nat1cup}}{\equref{qybe}}
\gbeg{9}{15}
\gvac{1}\got{1}{X}\gvac{2}\got{1}{Y}\gvac{2}\got{2}{A}\gnl
\gvac{1}\gcl{1}\gvac{2}\gcl{2}\gvac{2}\gcmu\gnl
\gvac{1}\gcl{1}\gvac{1}\gdb\gvac{2}\gcn{1}{1}{1}{-1}\gcl{1}\gnl
\gvac{1}\gcl{1}\gvac{1}\gcn{1}{1}{1}{0}\gcl{1}\glcm\glcm\gnl
\gvac{1}\gcl{1}\gcmu\gbr\gbr\gcl{1}\gnl
\gvac{1}\gcl{1}\gcl{1}\gcl{1}\gcl{1}\gbr\gcl{1}\gcl{1}\gnl
\gvac{1}\gcl{1}\gcl{1}\gcl{1}\gmu\gbr\gcl{1}\gnl
\gvac{1}\gcl{1}\gcl{1}\gcl{1}\gcn{1}{1}{2}{1}\gvac{1}\gcl{1}\gcl{7}\gvac{-1}\gev\gnl
\gvac{1}\gcl{1}\gcl{1}\gbr\gvac{1}\gcn{1}{1}{1}{-1}\gnl
\gvac{1}\gcl{1}\gmu\glm\gnl
\gdb\gvac{-1}\gcl{1}\gcn{1}{1}{2}{1}\gvac{1}\gcn{1}{2}{3}{-1}\gnl
\gcl{1}\gbr\gnl
\gmu\gcl{2}\gvac{-1}\gev\gnl
\gcn{1}{1}{2}{2}\gnl
\gob{2}{B}\gob{1}{X}\gvac{4}\gob{1}{Y}
\gend\\
&&\equalupdown{\equref{nat1cup}}{\equref{lcomcoalg}}
\gbeg{8}{12}
\gvac{1}\got{1}{X}\gvac{2}\got{1}{Y}\gvac{1}\got{1}{A}\gnl
\gvac{1}\gcl{1}\gvac{2}\gcl{2}\glcm\gnl
\gvac{1}\gcl{1}\gvac{1}\gdb\gcl{1}\gcn{1}{1}{1}{2}\gnl
\gvac{1}\gcl{1}\gvac{1}\gcn{1}{1}{1}{0}\gcl{1}\gcl{1}\gcmu\gnl
\gvac{1}\gcl{1}\gcmu\gbr\gcl{1}\gcl{1}\gnl
\gvac{1}\gcl{1}\gcl{1}\gbr\gbr\gcl{1}\gnl
\gvac{1}\gcl{1}\gmu\glm\gcl{5}\gvac{-1}\gev\gnl
\gdb\gvac{-1}\gcl{1}\gcn{1}{1}{2}{1}\gvac{1}\gcn{1}{2}{3}{-1}\gnl
\gcl{1}\gbr\gnl
\gmu\gcl{2}\gvac{-1}\gev\gnl
\gcn{1}{1}{2}{2}\gnl
\gob{2}{B}\gob{1}{X}\gvac{3}\gob{1}{Y}
\gend
\equalupdown{\equref{nat1cup}}{\equref{nat2cup}}
\gbeg{8}{13}
\gvac{1}\got{1}{X}\gvac{1}\got{1}{Y}\gvac{2}\got{1}{A}\gnl
\gvac{1}\gcl{2}\gvac{1}\gcl{2}\gvac{1}\glcm\gnl
\gdb\gdb\gvac{1}\gcl{1}\gcl{1}\gnl
\gcl{1}\gbr\gcn{1}{1}{1}{3}\gvac{1}\gcl{1}\gcl{1}\gnl
\gmu\gcl{1}\gvac{1}\gbr\gcn{1}{1}{1}{2}\gnl
\gcn{1}{2}{2}{3}\gvac{1}\gcl{1}\gvac{1}\gcn{1}{1}{1}{-1}\gcl{2}\gcmu\gnl
\gvac{2}\gbr\gdb\gcl{1}\gcl{1}\gnl
\gvac{1}\gmu\gcl{1}\gcl{1}\gbr\gcl{1}\gnl
\gvac{1}\gcn{1}{4}{2}{2}\gvac{1}\gcl{1}\glm\gcl{4}\gvac{-1}\gev\gnl
\gvac{3}\gcl{1}\gvac{1}\gcn{1}{1}{1}{-1}\gnl
\gvac{3}\gcl{2}\gvac{-1}\gev\gnl
\gvac{1}\gob{2}{B}\gob{1}{X}\gvac{2}\gob{1}{Y}
\gend
\hspace{1mm},
\end{eqnarray*}
and this shows that \equref{tensprodcomprel} holds. In this computation, we freely
used the fact that $Y$ is left $B$-comodule, that the multiplication on $B$ is associative, 
etc. \hfill$\square$

\section{Examples}\selabel{3}
\setcounter{equation}{0}
Let $\Cc={}_k{\cal M}$, the category of vector spaces over a field $k$. 
A bialgebra in ${}_k{\cal M}$ is an ordinary $k$-bialgebra, and, for a left $B$-comodule algebra 
$A$, the category ${}^B\Cc_A={}^B{\cal M}_A$ is the classical category of relative $(B, A)$ Hopf modules, see \cite{mont, tak}. 
It is well-known that particular examples of 
bialgebras (even Hopf algebras) in ${}^B_B{\cal YD}$ can be obtained from quasitriangular or 
coquasitriangular Hopf algebras. We will study these two classes of examples in more detail.\\
For the definitions of a quasitriangular and a coquasitriangular Hopf algebra, we 
invite the reader to consult \cite{kas, mont, maj}. In the sequel, we denote the $R$-matrix of a 
quasitriangular Hopf algebra 
$H$ by $R=R^1\ot R^2\in H\ot H$, and the bilinear form that defines a 
coquasitriangular structure on a Hopf algebra $H$ by 
$\sigma :\ H\ot H\ra k$.\\
Let $(H, R)$ be a quasitriangular Hopf algebra with antipode $S$. The enveloping algebra
braided group $\un{H}$ is equal to $H$ as an algebra, has the same unit and counit,
but newly defined comultiplication $\un{\Delta}$:
\begin{equation}\eqlabel{comEABG}
\un{\Delta}(h)=h_1S(R^2)\ot R^1\tr h_2,
\end{equation} 
where $\tr$ is the left adjoint action, that is, $h\tr h'=h_1h'S(h_2)$, for all $h, h'\in H$. 
Now $\un{H}$ is a braided bialgebra in ${}_H{\cal M}$, the category of left representations
of $H$, see \cite[Ex. 9.4.9]{maj}.
Now $(H, R)$ is quasitriangular, so we have a braided functor
$F: {}_H{\cal M}\ra {}^H_H{\cal YD}$, see \cite[Lemma 7.4.4]{maj}. For $M\in {}_M{\cal M}$,
$F(M)=M$ with its original left $H$-action, and left $H$-coaction defined by
\begin{equation}\eqlabel{coactEABG}
\l_M: M\ra H\ot M~,~\l_M(m)=R^2\ot R^1\cdot m~,~m\in M~.
\end{equation}
A braided functor sends bialgebras to bialgebras, hence it follows that
$\un{H}$ is a bialgebra in ${}^H_H{\cal YD}$. 
This provides the following example.

\begin{example}
Let $(H, R)$ be a quasitriangular Hopf algebra with antipode $S$ and $\un{H}$ the enveloping algebra braided group of 
$H$. Then the category of relative Hopf modules ${}^H{\cal M}_{\un{H}}$ is a monoidal category. The tensor product 
is the usual tensor product of vector spaces endowed with the left $H$-comodule structure 
given by the comultiplication $\Delta$ of $H$ and with the right $\un{H}$-module structure given by 
\[
(x\ot y)\bullet h=x\cdot (y_{-1}\tr h_1S(R^2))\ot y_0\cdot (R^1\tr h_2)~,
\]
for all $x\in X\in {}^H{\cal M}_{\un{H}}$ and $y\in Y\in {}^H{\cal M}_{\un{H}}$ (here $Y\ni y\mapsto y_{-1}\ot y_0\in H\ot Y$ 
is the Sweedler notation for the left $H$-coaction on $Y$). The unit object is $k$ considered as a left $H$-comodule via 
the unit of $H$ and a right $\un{H}$-module via the counit of $H$, and the associativity, left and right unit 
constraints are those of ${}_k{\cal M}$.  
\end{example}

\begin{proof}
We only point out that an object $X$ of ${}^H{\cal M}_{\un{H}}$ is a left $H$-comodule and a right $H$-module ($\un{H}=H$ as 
algebras) for which the following compatibility relation holds:
\[
(x\cdot h)_{-1}\ot (x\cdot h)_0=x_{-1}R^2\ot x_0\cdot (R^1\tr h)~,
\]
for all $x\in X$ and $h\in H$.
\end{proof}

\begin{corollary}
If $H$ is a cocommutative Hopf algebra, then the category of right-left Long $H$-dimodules is a monoidal category.
\end{corollary}

\begin{proof}
If $H$ is cocommutative, then it is quasitriangular with $R=1\ot 1$. In this particular situation, we have
that $\un{H}=H$ as an ordinary bialgebra, 
and an object $X$ of ${}^H{\cal M}_{\un{H}}$ is a left $H$-comodule and a right $H$-module 
such that 
\[
(x\cdot h)_{-1}\ot (x\cdot h)_0=x_{-1}\ot x_0\cdot h,
\]
for all $x\in X$ and $h\in H$.
Comparing to the right-left version of \cite[Def. 16]{stefmz} we can conclude that, in this situation, 
the category of right-left $H$-dimodules is ${}^H{\cal M}_{\un{H}}$, and so it is a monoidal category. 
Observe that the right $H$-module structure on the tensor product 
is given by the formula
\[
(x\ot y)\bullet h=x\cdot (y_{-1}\tr h_1)\ot y_0\cdot h_2~,
\]
for any two right-left $H$-dimodules $X,Y$, and any $x\in X$, $y\in Y$ and $h\in H$. 
\end{proof}

We now move to the dual situation: let $(H, \sigma)$ be a coquasitriangular Hopf algebra with antipode $S$ 
($S$ is then bijective, see \cite{doi}). The (left) function algebra braided group $\un{H}$
is equal to $H$ as a coalgebra, with the same unit and counit, and multiplication 
$\diamond$ defined by 
\begin{equation}\eqlabel{multFABG}
h\diamond h'=\sigma(h'_2, S(h_1)h_3)h_2h'_1,
\end{equation}
for all $h, h'\in H$. Note that we mention explicitely that we consider the left handed version
of the function algebra braided group; it is obtained from $H$ using the left version
of the transmutation theory, see \cite[Remark 4.3]{bb}. In \cite[Ex. 9.4.10]{maj}, the
right handed version is presented. We need the left version here since it fits in our
context; namely $\un{H}$ is a braided bialgebra in ${}^H{\cal M}$ via the left adjoint
coaction $\lambda$ given by the formula
$\lambda(h)=S^{-1}(h_3)h_1\ot h_2$, for all $h\in H$. The multiplication is given by 
\equref{multFABG}, and the other structure maps on $\un{H}$ coincide with the corresponding ones
on $H$.\\
$(H, \sigma)$ is a coquasitriangular Hopf algebra, hence there exists a braided monoidal
functor $G: {}^H{\cal M}\ra {}^H_H{\cal YD}$. At the level of objects, $G(M)=M$,
with the original left $H$-coaction, and left $H$-action given by the formula
$h\cdot m= \sigma(m_{-1}, h)m_{0}$, for all $h\in H$ and $m\in M$. 
Consequently $\un{H}$ is a bialgebra in ${}^H_H{\cal YD}$ 
via the left coadjoint coaction and $H$-action defined by $h\succ h'=\sigma(S^{-1}(h'_3)h'_1, h)h'_2$, for all $h, h'\in H$. 
We then obtain the following result.

\begin{example}
Let $(H, \sigma)$ be a coquasitriangular Hopf algebra and $\un{H}$ the 
associated (left) function algebra braided group. 
Then ${}^H{\cal M}_{\un{H}}$ is a monoidal category with tensor product inherited from ${}^H{\cal M}$ and 
equipped with the additional right $\un{H}$-module structure given by 
\[
(x\ot y)\bullet h=\sigma(S^{-1}(h_3)h_2, y_{-1})x\cdot h_1\ot y_0\cdot h_4~,
\] 
for all $x\in X\in {}^H{\cal M}_{\un{H}}$, $y\in Y\in {}^H{\cal M}_{\un{H}}$ and $h\in \un{H}$. 
The unit object is $k$ considered as a left $H$-comodule via 
the unit of $H$ and a right $\un{H}$-module via the counit of $H$, and the associativity, left and right unit constraints are those of ${}_k{\cal M}$.  
\end{example}

\begin{proof}
Everything follows from the above considerations and results. Observe that an object in 
${}^H{\cal M}_{\un{H}}$ is an ordinary left $H$-comodule $M$ which is at the same time a right $\un{H}$-module, that is, 
$m\cdot 1=m$ and 
\[
(m\cdot h)\cdot h'=\sigma(h'_2, S(h_1)h_3)m\cdot (h_2h'_1)~, 
\]
for all $m\in M$, $h, h'\in H$, and such that the following compatibility relation holds,
\[
(m\cdot h)_{-1}\ot (m\cdot h)_0=m_{-1}S^{-1}(h_3)h_1\ot m_0\cdot h_2~,
\]  
for all $m\in M$ and $h\in H$. 
\end{proof}

\begin{remark}
A commutative Hopf algebra $H$ is coquasitriangular with trivial $\sigma$, $\sigma(h, h')=\va(h)\va(h')$. In this case, the 
(left) function algebra braided group associated to $H$ is $H$ itself 
with the left coadjoint coaction. Therefore, an object $M$ in ${}^H{\cal M}_{\un{H}}$ is a left $H$-comodule and a right $H$-module 
such that 
\[
(m\cdot h)_{-1}\ot (m\cdot h)_0=m_{-1}S^{-1}(h_3)h_1\ot m_0\cdot h_2~,
\]  
for all $m\in M$ and $h\in H$. Since $H$ is commutative it follows that the above condition 
is equivalent to the required compatibility relation for a right-left Yetter-Drinfeld module over $H$. Hence ${}^H{\cal YD}_H={}^H{\cal M}_{\un{H}}$. 
It is well-known, of course, that the category of Yetter-Drinfeld modules is monoidal;
but perhaps it is interesting to know that the category of Yetter-Drinfeld modules over a
commutative Hopf algebra can be identified with a suitable category of relative
Hopf modules.
\end{remark}

Finally, we discuss the relationship with the monoidal structures that were discussed in \cite{cob}.
Let $A$ be a left $B$-module algebra. Then $(B,A,B)$ is a Doi-Hopf datum in the sense of
\cite{doi2}. Monoidal Doi-Hopf data were introduced in \cite{cob}, and it is easy to see that
$(B,A,B)$ is monoidal if and only if $A$ is a bialgebra (in the category of vector spaces), and
\begin{eqnarray}
ha_{[-1]}\ot \Delta(a_{[0]})&=& a_{(1)[-1]}ha_{(2)[-1]}\ot a_{(1)[0]}\ot a_{(2)[0]}\eqlabel{cob1};\\
\varepsilon_A(a)1_B&=& \varepsilon_A(a_{[0]})a_{[-1]}\eqlabel{cob2},
\end{eqnarray}
for all $h\in B$ and $a\in A$.
Note that the left-right 
convention in \cite{cob} is different from ours.
We have used the Sweedler notation, where indices between brackets refer to comultiplication
and indices between square brackets refer to coaction.

\begin{proposition}\prlabel{3.4}
Let $A$ be a left $B$-comodule algebra, equipped with the trivial left $B$-action
$h\cdot a=\varepsilon(h)a$. Then $(B,A,B)$ is a monoidal Doi-Hopf datum if and only
if $A$ is bialgebra in ${}_B^B{\cal YD}$.
\end{proposition}

\begin{proof}
$A$ is a braided algebra if and only
the following conditions hold,
\begin{enumerate}
\item $A$ is a left Yetter-Drinfeld module over $B$;
\item $A$ is a left $B$-(co)module (co)algebra;
\item $\varepsilon_A$ is an algebra map and $\Delta_A(1_A)=1_A\ot 1_A$;
\item $\Delta_A(ab)=a_{(1)}(a_{(2)[-1]}\cdot b_{(1)}) \ot a_{(2)[0]}b_{(2)}$.
\end{enumerate}
First assume that $(B,A,B)$ is monoidal. (3) is satisfied since $A$ is a bialgebra;
(4) simplifies to $\Delta(ab)=\Delta(a)\Delta(b)$, since the $B$-action on $A$ is trivial,
and this is also satisfied. Three of the four conditions in (2) are satisfied, the only one that is
left to prove is the fact that $A$ is a left $B$-comodule coalgebra. This follows from \equref{cob2}
and \equref{cob1} (with $h=1_B$). Applying $\varepsilon_A$ to \equref{cob1}, we find
$$ha_{[-1]}\ot a_{[0]}=a_{(1)[-1]}ha_{(2)[-1]}\varepsilon_A(a_{(2)[0]})\ot a_{(1)[0]}
\equal{\equref{cob2}} a_{[-1]}h\ot a_{[0]},$$
which is precisely the compatibility relation for Yetter-Drinfeld modules, at least in the case
where the $H$-action is trivial.\\
Conversely, assume that $A$ is a braided bialgebra. $A$ is a left $B$-comodule coalgebra,
so \equref{cob2} holds and
\begin{equation}\eqlabel{cob5}
a_{[-1]}\ot \Delta(a_{[0]})= a_{(1)[-1]}a_{(2)[-1]}\ot a_{(1)[0]}\ot a_{(2)[0]}.
\end{equation}
Furthermore,
\begin{equation}\eqlabel{cob6}
ha_{[-1]}\ot a_{[0]}=a_{[-1]}h\ot a_{[0]},
\end{equation}
since $A$ is a Yetter-Drinfeld module. Then
\begin{eqnarray*}
&&\hspace*{-2cm}
ha_{[-1]}\ot \Delta(a_{[0]})\equal{\equref{cob6}}
a_{[-1]}h\ot \Delta(a_{[0]})
\equal{\equref{cob5}}
a_{(1)[-1]}a_{(2)[-1]}h\ot a_{(1)[0]}\ot a_{(2)[0]}\\
&\equal{\equref{cob6}}&
a_{(1)[-1]}ha_{(2)[-1]}\ot a_{(1)[0]}\ot a_{(2)[0]}.
\end{eqnarray*}
This proves that \equref{cob1} holds, and the result follows.
\end{proof}

If $(B,A,B)$ is a monoidal Doi-Hopf data, then we have a monoidal structure on the category
of relative Doi-Hopf modules, see \cite[Prop. 2.1]{cob}. This monoidal structure coincides
with the one that follows from \thref{2.1}.

\begin{example}
We end with a trivial example. Let $B=k$, and $A$ a $k$-algebra. The category of Yetter-Drinfeld
modules ${}_k^k{\cal YD}$ is just the category of vector spaces ${\cal M}_k$, and a bialgebra in this
category is an ordinary bialgebra. So we recover the classical result that monoidal structures
on the category of representations of an algebra $A$ are in one-to-one correspondence
with bialgebra structures on $A$.
\end{example}


\end{document}